\documentclass[10pt]{article}
\usepackage{amsfonts}
\usepackage{amsmath,amssymb}
\usepackage[a4paper]{geometry}

\newtheorem{Lemma}{Lemma}[section]
\newtheorem{theorem}[Lemma]{Theorem}

\newtheorem{lemma}[Lemma]{Lemma}

\newtheorem{definition}[Lemma]{Definition}

\begin{document}

\begin{center}
{\Large \textbf{Well-posedness of stochastic second grade fluids}} \vspace{%
3mm}\\[0pt]

\vspace{1 cm} \textsc{\ Nikolai Chemetov}{\footnote{%
Departamento de Matem\'atica, Faculdade de Ci\^encias da Universidade de
Lisboa. E-mail: nvchemetov{\char'100}fc.ul.pt, nvchemetov{\char'100}%
gmail.com.},} \hspace{0mm} \textsc{Fernanda Cipriano}{\footnote{%
Departamento de Matem\'atica, Faculdade de Ci\^encias e Tecnologia da
Universidade Nova de Lisboa and Centro de Matem\'atica e Aplica\c c\~oes.
E-mail: cipriano{\char'100}fct.unl.pt.}}


\end{center}

\noindent 

\date{\today }

\begin{abstract}
The theory of turbulent Newtonian fluids turns out that the choice of the
boundary condition is a relevant issue, since it can modify the behavior of
the fluid by creating or avoiding a strong boundary layer. In this work we
study stochastic second grade fluids filling a two-dimensional bounded
domain, with the Navier-slip boundary condition (with friction). We prove
the well-posedness of this problem and establish a stability result. Our
stochastic model involves a multiplicative white noise and a convective term
with third order derivatives, which significantly complicate the analysis.
\vspace{2mm} \newline

\vspace{3mm}\newline


\vspace{3mm}
\end{abstract}

\section{Introduction}

\setcounter{equation}{0}
The present work is devoted to the study of the stochastic incompressible
fluids of second grade, which are a special class of non-Newtonian fluids.
Unlike the Newtonian fluids, where only the stretching tensor appears in the
characterization of the stress response to a deformation fluid, here the
Cauchy stress tensor $\mathbb{T}$ of the non-Newtonian fluids is defined by
\begin{equation*}
\mathbb{T}=-\pi \mathbb{I}+\nu A_{1}+\alpha _{1}A_{2}+\alpha _{2}\,A_{1}^{2},
\end{equation*}%
where the first term $-\pi \mathbb{I}$ is due to the incompressibility of
the fluid and $\mathbb{A}_{1}$, $\mathbb{A}_{2}$ are the two first
Rivlin-Ericksen tensors (cf. \cite{RE55})
\begin{equation*}
\mathbb{A}_{1}(y)=\nabla {y}+\left( \nabla {y}\right) ^{\top }\quad \text{%
and }\quad \mathbb{A}_{2}(y)=\dot{\mathbb{A}}_{1}(y)+\mathbb{A}_{1}(y)\nabla
y+\left( \nabla y\right) ^{\top }\mathbb{A}_{1}(y),
\end{equation*}%
where $y$ denotes the velocity of the fluid, the superposed dot is the
material time derivative, $\nu $ is the kinematic viscosity of the fluid and
$\alpha _{1}$, $\alpha _{2}$ are constant material moduli. The study
developed in \cite{DF74} turns out that thermodynamic laws and stability
principles impose $\alpha _{1}\geq 0$ and $\alpha _{1}+\alpha _{2}=0.$ We
set $\alpha=\alpha _{1}$ and assume $\alpha _{1}>0.$

It is well known that in turbulent fluids, small random perturbations can
produce relevant macroscopic effects. By this reason, the incorporation of a
stochastic white noise force in the Navier-Stokes equations \cite{BT73} is
widely recognized as an important step to understand the turbulence
phenomena. In this perspective, we can find in \cite{B95} (see Lemma 2.2) a
deduction of stochastic Navier-Stokes equations from fundamental principles,
by showing that the stochastic Navier-Stokes equations are a real physical
model. Nowadays, the stochastic Navier-Stokes equations are quite well
understood, see for instance in \cite{D08}, \cite{Fla08}, \cite{MR05}, \cite%
{S10} and the references therein. In spite of that, there are few results in
the literature about stochastic non-Newtonian fluids \cite{D11}, \cite{RS101}%
, \cite{RS102}, \cite{RS12}. In this paper we consider the stochastic second
grade equations with multiplicative noise given by
\begin{equation}
\left\{
\begin{array}{l}
\tfrac{\partial }{\partial t}\left( Y-\alpha \Delta Y\right) =\nu \Delta Y-%
\mathrm{curl}\left( Y-\alpha \Delta Y\right) \times Y-\nabla \pi +U+G(t,Y)\,%
\dot{W}_{t},\vspace{2mm} \\
\mathrm{div}\,Y=0\qquad \qquad \qquad \qquad \qquad \qquad \qquad \qquad
\qquad \qquad \qquad \qquad \qquad \text{in }\mathcal{O}\times (0,T),\vspace{%
2mm}%
\end{array}%
\right.  \label{equation_etat_temps}
\end{equation}%
where $U$ is a body force, $G(t,Y)\,\dot{W}_{t}$ is a multiplicative white
noise and $\mathcal{O}$ is a bounded domain of $\mathbb{R}^{2}$ with a
boundary $\Gamma $.

The study of this system requires suitable boundary conditions on the
boundary $\Gamma $ of the domain. The Dirichlet boundary condition given by
\begin{equation*}
\begin{array}{ll}
Y=0 & \qquad \mbox{on}\ \Gamma%
\end{array}%
\end{equation*}%
is accepted as an appropriate boundary condition and is the more usual one.
Another physical relevant boundary condition considered in the literature is
the Navier boundary condition
\begin{equation}
Y\cdot \mathrm{n}=0,\qquad \left[ 2(\mathrm{n}\cdot DY)+\gamma Y\right]
\cdot \mathrm{\tau }=0\qquad \mbox{on}\ \Gamma ,  \label{NS}
\end{equation}%
where $\mathrm{n}=(n_{1},n_{2})$ and \textrm{$\tau $}$=(-n_{2},n_{1})$ are
the unit normal and tangent vectors, respectively, to the boundary $\Gamma $%
, $DY=\tfrac{\nabla Y+(\nabla Y)^{\top }}{2}$ is the symmetric part of the
velocity gradient
and $\gamma >0$ is a friction coefficient on $\Gamma .$

The stochastic partial differential equations (\ref{equation_etat_temps}) with the
Dirichlet boundary condition has been studied
 in \cite{RS101} and \cite{RS12}. 
In the former paper, the authors used tightness arguments that conjugated
with the Skorohod theorem provided the existence of a weak stochastic
solution, in the sense that the Brownian motion, being part of the solution,
was not given in advance; while in the second one, the authors proved the
existence and uniqueness of a strong stochastic solution. Let us refer the
pioneer papers \cite{O81} and \cite{CO84} (see also \cite{CG97}), where the
deterministic second grade equations with the Dirichlet boundary condition
were mathematically studied for the first time, and \cite{BR03} where the
deterministic equations were studied with a particular Navier boundary
condition (without friction, i.e. when $\gamma =0$). The physical
interpretation of these second grade equations can be found in \cite{CH93},
\cite{DF74}, \cite{DR95}, \cite{FHT2}, \cite{HMR981} and \cite{HMR98}. It is
relevant to recall that the deterministic methods are based on the
Faedo-Galerkin approximation method and a priori estimates. Then,
compactness arguments can be used to pass to the limit of the respective
approximate equations in the distributional sense. Unfortunately, for the
stochastic partial differential equations a priori estimates are not enough to pass
to the limit of the approximate equations, due to the lack of regularity on
the time and stochastic variables. In order to obtain a strong stochastic
solution we should verify that the sequence of the Galerkin approximations
converges strongly in some adequate topology.

We should mention that even if the Dirichlet boundary condition is widely
accepted as an appropriate boundary condition at the surface of contact
between a fluid and a solid, it is also a source of many problems since it
attaches fluid particles to the boundary, creating a strong boundary layer
(cf. \cite{CMR98}, \cite{IP06}, \cite{K06}, \cite{LNP05}). On the other
hand, the Navier boundary condition allows the slippage of the fluid on the
boundary, making it possible to treat important problems as for instance the
boundary layer problem, when the viscosity $\nu $ and/or the elastic
response $\alpha $ tend to zero (cf. \cite{BILN12}, \cite{cc1}-\cite{cc3},
\cite{CT15}, \cite{LT10}, \cite{LNTZ15}). However, even if the Navier-slip
boundary condition allows to solve interesting problems, technically, when
comparing with the Dirichlet boundary condition, it requires a more careful
mathematical analysis to show the well-posedness of system (\ref%
{equation_etat_temps})-(\ref{NS}) as well as to establish stability
properties for the solution, since the boundary terms resulting from
integrating by parts of the convective term do not vanish and should be
estimated in an appropriate way.

As far as we know, the stochastic second grade fluid equations with the
Navier boundary condition are studied here for the first time. To show the
well-posedness, as in previous articles, we follow the Faedo-Galerkin
approximation method by taking an appropriate basis. We first deduce uniform
estimates for the approximate solutions that allow to pass to the limit with
respect to the weak topology. In order to show that the limit process is a
solution, we adapt the methods developed in \cite{B99} to study the
stochastic Navier-Stokes equations. More precisely, we show that the
approximate solutions already converge strongly up to a certain stopping
time, therefore we establish the existence and uniqueness results for the
solution of system (\ref{equation_etat_temps})-(\ref{NS}), as a stochastic
process with values in $H^{3}$. We should mention that an analogous
reasoning is considered in \cite{RS12} to deal with the stochastic second
grade fluid equations with homogeneous Dirichlet boundary condition.

The plan of the present paper is as follows. In Section 2 we state the
functional setting and introduce useful notations. In Section 3 we present
some well known results and relevant lemmas related with the nonlinear term
of (\ref{equation_etat_temps})$_{1}$, which will be applied in the next
sections. The main result concerning the existence of a strong stochastic
solution is established in Section 4. Finally Section 5 is devoted to the
study of the stability property.



\section{Functional setting and notations}

\setcounter{equation}{0}
We consider the stochastic second grade fluid model in a bounded and simply
connected domain $\mathcal{O}$ of $\mathbb{R}^{2}$ with a sufficiently
regular boundary $\Gamma $
\begin{equation}
\left\{
\begin{array}{cc}
d(\upsilon (Y))=(\nu \Delta Y-\mathrm{curl}(\upsilon (Y))\times Y-\nabla \pi
+U)\,dt+G(t,Y)\,dW_{t}, &  \\
\multicolumn{1}{l}{\mathrm{div}\,Y=0} & \multicolumn{1}{l}{\mbox{in}\
\mathcal{O}\times (0,T),} \\
\multicolumn{1}{l}{Y\cdot \mathrm{n}=0,\qquad \left[ 2(\mathrm{n}\cdot
DY)+\gamma Y\right] \cdot \mathrm{\tau }=0} & \multicolumn{1}{l}{\mbox{on}\
\Gamma \times (0,T),} \\
\multicolumn{1}{l}{Y(0)=Y_{0}} & \multicolumn{1}{l}{\mbox{in}\ \mathcal{O},}%
\end{array}%
\right.  \label{equation_etat}
\end{equation}%
where $\nu >0$ is a constant viscosity of the fluid, $\alpha >0$ is a
constant material modulus, the constant $\gamma >0$ is a friction
coefficient of $\Gamma ,$ $\ \Delta $ and $\nabla $ respectively denote the Laplacian and
the gradient, $Y=(Y_{1},Y_{2})$ is a 2D velocity field and
\begin{equation*}
\upsilon (Y)=Y-\alpha \Delta Y.
\end{equation*}%
The function $\pi $ represents the pressure, $U$ is a distributed mechanical
force and the term 
$$G(t,Y) \,dW_{t}=\sum_{k=1}^m  G^k(t,Y) \,dW^k_{t}$$ 
corresponds to the stochastic
perturbation, where $G(t,Y)=(G^1(t,Y), \dots, G^m (t,Y))$ has 
suitable growth assumptions defined below and  $W_{t}=(W^1_{t}, \dots, W^m_{t})$ is a standard $\mathbb{R}^{m}$-valued Wiener
process defined on a complete probability space $(\Omega ,\mathcal{F},P)$
endowed with a filtration $\left\{ \mathcal{F}_{t}\right\} _{t\in \lbrack
0,T]}$. We assume that $\mathcal{F}_{0}$\ contains every $P$-null subset of $%
\Omega $.

Let $X$ be a real Banach space endowed with the norm $\left\Vert \cdot
\right\Vert _{X}.$ We denote by $L^{p}(0,T;X)$ the space of $X$-valued
measurable $p-$integrable functions $y$ defined on $[0,T]$ for $p\geq 1$.

For $p,r\geq 1$ let $L^{p}(\Omega ,L^{r}(0,T;X))$ be the space of processes $%
y=y(\omega ,t)$ with values in $X$ defined on $\ \Omega \times \lbrack 0,T],$
adapted to the filtration $\left\{ \mathcal{F}_{t}\right\} _{t\in \lbrack
0,T]}$ ,\ and endowed with the norms
\begin{equation*}
\left\Vert y\right\Vert _{L^{p}(\Omega ,L^{r}(0,T;X))}=\left( \mathbb{E}%
\left( \int_{0}^{T}\left\Vert y\right\Vert _{X}^{r}\,dt\right) ^{\frac{p}{r}%
}\right) ^{\frac{1}{p}}\text{ }
\end{equation*}%
and%
\begin{equation*}
\left\Vert y\right\Vert _{L^{p}(\Omega ,L^{\infty }(0,T;X))}=\left( \mathbb{E%
}\sup_{t\in \lbrack 0,T]}\left\Vert y\right\Vert _{X}^{p}\ \right) ^{\frac{1%
}{p}}\quad \text{if }r=\infty ,
\end{equation*}%
where $\mathbb{E}$ is the mathematical expectation with respect to the
probability measure $P.$ As usual in the notation of processes $y=y(\omega
,t)$ we normally omit the dependence on $\omega \in \Omega .$

In equation (\ref{equation_etat}) the vector product $\times $ for 2D
vectors $y=(y_{1},y_{2})$ and $z=(z_{1},z_{2})$ \ is calculated as $y\times
z=(y_{1},y_{2},0)\times (z_{1},z_{2},0);$ the curl of the vector $y$ is
equal to $\mathrm{curl}\,y=\tfrac{\partial y_{2}}{\partial x_{1}}-\tfrac{%
\partial y_{1}}{\partial x_{2}}$ and the vector product of $\mathrm{curl}\,y$
with the vector $z$ is understood as
\begin{equation*}
\mathrm{curl}\,y\times z=(0,0,\mathrm{curl}\,y)\times (z_{1},z_{2},0).
\end{equation*}%
Given two vectors $y,z\in \mathbb{R}^{2}$, $y\cdot
z=\sum_{i=1}^{2}y_{i}z_{i} $ stands for the usual scalar product in $\mathbb{%
R}^{2}$ and given two matrices $A, B,$ we denote $A\cdot B=\sum_{ij=1}^2
A_{ij}B_{ij}$.

Let us introduce the following Hilbert spaces
\begin{equation}
\begin{array}{l}
H(\mathrm{curl};\mathcal{O})=\left\{ y\in L^{2}(\mathcal{O})\mid \mathrm{curl%
}\,y\in L^{2}(\mathcal{O}),\quad \mathrm{div}\,y=0\ \text{ in }\mathcal{O}%
\right\} ,\vspace{2mm} \\
H=\left\{ y\in L^{2}(\mathcal{O})\mid \mathrm{div}\,y=0\ \text{ in }\mathcal{%
O}\ \mbox{ and }\ y\cdot \mathrm{n}=0\ \mbox{ on }\Gamma \right\} ,\vspace{%
2mm} \\
V=\left\{ y\in H^{1}(\mathcal{O})\mid \mathrm{div}\,y=0\ \mbox{ in }\
\mathcal{O}\mbox{ and }\ y\cdot \mathrm{n}=0\ \text{ on }\ \Gamma \right\} ,%
\vspace{2mm} \\
W=\left\{ y\in V\cap H^{2}(\mathcal{O})\mid \left[ 2(\mathrm{n}\cdot
Dy)+\gamma y\right] \cdot \mathrm{\tau }=0\ \ \mbox{on}\ \Gamma \right\} ,%
\vspace{2mm} \\
\widetilde{W}=W\cap H^{3}(\mathcal{O}).%
\end{array}
\label{w}
\end{equation}%
We denote by $(\cdot ,\cdot )$ the inner product in $L^{2}(\mathcal{O})$ and
by $\Vert \cdot \Vert _{2}$ the associated norm. The norm in the space $%
H^{p}(\mathcal{O})$ is denoted by $\Vert \cdot \Vert _{H^{p}}$. Let us note
that $H(\mathrm{curl};\mathcal{O})$ is a subspace of $H^{1}(\mathcal{O}).$
Let us denote
\begin{equation*}
\left( Dy,Dz\right) =\int_{\mathcal{O}}Dy\cdot Dz.
\end{equation*}%
On the space $V$, we consider the following inner product
\begin{equation*}
\left( y,z\right) _{V}=\left( \upsilon (y),z\right) =\left( y,z\right)
+2\alpha \left( Dy,Dz\right) +\alpha \gamma \int_{\Gamma }y\cdot z
\end{equation*}%
and the corresponding norm $\Vert \cdot \Vert _{V}.$ We can verify that the
norms $\Vert \cdot \Vert _{H^{1}}$ and $\Vert \cdot \Vert _{V}$ are
equivalent because of the Korn inequality
\begin{equation}
\left\Vert y\right\Vert _{H^{1}}\leq C\left( \left\Vert Dy\right\Vert
_{2}+\left\Vert y\right\Vert _{2}\right) ,\qquad \forall y\in H^{1}(\mathcal{%
O}).  \label{korn}
\end{equation}%
Here and below, $C$ will denote a generic positive constant that may depend
only on the domain $\mathcal{O},$ \ the regularity of the boundary $\Gamma $%
, the physical constants $\nu $, $\alpha $, $\gamma $ and $K$, defined in \ (%
\ref{LG}).

Let $B$ be a given Hilbert space  with inner product $(\cdot.\cdot)_B$. For a vector 
$$
h=(h^1,\dots,h^m)\in B^m= \overbrace{B\times ...\times B}^{m-times}
$$ 
we introduce the norm 
\begin{equation*}
\left\Vert h\right\Vert _{B}=\sum_{i=1}^{m}\left\Vert
h_{i}\right\Vert _{B} 
\end{equation*}
and the module of the inner product of $h$ and a fixed $v\in B$ as
\begin{equation}
|\left(h,v\right)_B |  = \left( \sum_{k=1}^{m}\left( h^k, v \right)_B^2 \right)^{1/2}.  \label{product}
\end{equation}. 

Assume that $G(t,y): [0,T]\times V\rightarrow  V^{m} $ is Lipschitz on $y$, and satisfies a linear growth; that is, there exists a
 positive constant $K$ such that
\begin{align}
\left\Vert G(t,y)-G(t,z)\right\Vert _{V}^{2}& \leq K\left\Vert
y-z\right\Vert _{V}^{2},  \notag \\
\left\Vert G(t,y)\right\Vert _{V}& \leq K\left( 1+\left\Vert y\right\Vert
_{V}\right) ,\qquad \forall y,z\in V,\;t\in \lbrack 0,T].  \label{LG}
\end{align}%

\section{ Preliminary results}

\bigskip

\setcounter{equation}{0}
Let us introduce the Helmholtz projector $\mathbb{P}:L^{2}(\mathcal{O}%
)\longrightarrow H$, which is the linear bounded operator defined by $%
\mathbb{P}y=\tilde{y}$, where $\tilde{y}\in H$ is characterized by the
Helmholtz decomposition
\begin{equation*}
y=\tilde{y}+\nabla \phi ,\qquad \phi \in H^{1}(\mathcal{O}).
\end{equation*}%
We recall some useful inequalities, namely, the Poincar\'{e} inequality
\begin{equation*}
\left\Vert y\right\Vert _{2}\leq C\left\Vert \nabla y\right\Vert _{2}\qquad %
\mbox{for all}\ y\in V
\end{equation*}%
and the Sobolev inequality
\begin{equation*}
\left\Vert y\right\Vert _{4}\leq C\left\Vert \nabla y\right\Vert _{2}\qquad %
\mbox{for all}\ y\in V.
\end{equation*}

Now, we present the first result of this section. This is a well known and
very important property concerning the Navier boundary conditions (see Lemma
4.1 and Corollary 4.2 in \cite{K06}). \ Let $k$ be the curvature of $\Gamma $%
. Parameterizing $\Gamma $ by arc length $s$, the following relation holds
\begin{equation*}
\tfrac{\partial \mathrm{n}}{\partial \mathrm{\tau }}=\tfrac{d\mathrm{n}}{ds}%
=k\mathrm{\tau }.
\end{equation*}

\begin{lemma}
\label{curl}Let $y\in H^{2}(\mathcal{O})\cap V$ be a vector field verifying
the Navier boundary condition. Then
\begin{equation}
\mathrm{curl}\,y=g(y)\quad \text{ on}\ \Gamma \quad \text{ with }\
g(y)=\left( 2k-\gamma \right) \,y\cdot \mathrm{\tau }.  \label{curl_trace}
\end{equation}
\end{lemma}

\textbf{Proof.} Let us first notice that the anti-symmetric tensor $%
Ay=\nabla y-(\nabla y)^\top$ can be written in the form
\begin{equation*}
Ay=\mathrm{curl}\,y{\small \left[
\begin{array}{lr}
0 & -1 \\
1 & 0%
\end{array}%
\right] .}
\end{equation*}%
The symmetry of $Dy$ and the anti-symmetry of $Ay$ imply that
\begin{equation*}
\left( Dy\right) \mathrm{\tau }\cdot \mathrm{n}=\left( Dy\right) \mathrm{n}%
\cdot \mathrm{\tau }\quad \text{and}\quad \left( Ay\right) \mathrm{\tau }%
\cdot \mathrm{n}=-\left( Ay\right) \mathrm{n}\cdot \mathrm{\tau }.
\end{equation*}%
It follows that
\begin{equation*}
(\nabla y)\mathrm{\tau }\cdot \mathrm{n}=\left( Dy\right) \mathrm{n}\cdot
\mathrm{\tau }-\tfrac{1}{2}\left( Ay\right) \mathrm{n}\cdot \mathrm{\tau }
\end{equation*}%
which is equivalent to
\begin{equation}
\mathrm{curl}\,y=-2(\nabla y)\mathrm{\tau }\cdot \mathrm{n}+2(Dy)\mathrm{n}%
\cdot \mathrm{\tau }.  \label{navier}
\end{equation}%
Taking the derivative of the expression $y\cdot \mathrm{n}=0$ in the
direction of the tangent vector $\mathrm{\tau }$, we deduce
\begin{equation}
(\nabla y)\mathrm{\tau }\cdot \mathrm{n}=-k\,y\cdot \mathrm{\tau }.
\label{flux}
\end{equation}%
The conclusion is then a consequence of (\ref{navier}) and (\ref{flux}). $%
\hfill \hfill \blacksquare \hfill $

\bigskip

Now, we state a formula that can be easily derived by taking integration by
parts
\begin{equation}
-\int_{\mathcal{O}}\triangle y\cdot z=-\int_{\Gamma }2(Dy)n\cdot z+\int_{%
\mathcal{O}}2\,Dy\cdot Dz,  \label{integrate}
\end{equation}%
that holds for any $y\in H^{2}(\mathcal{O})\cap V$ and $z\in H^{1}(\mathcal{O%
})$. Using the boundary conditions, that gives the relation
\begin{equation}
-\int_{\mathcal{O}}\triangle y\cdot z=\gamma \int_{\Gamma }y\cdot z+\int_{%
\mathcal{O}}2\,Dy\cdot Dz\quad \text{ for any }y\in W\text{ and }z\in V
\label{integrate2}
\end{equation}%
that will be used throughout the article.

\bigskip

Let us consider the following modified Stokes system with Navier boundary
condition
\begin{equation}
\left\{
\begin{array}{ll}
h-\alpha \Delta h+\nabla p=f,\qquad \mathrm{div}\,h=0 & \quad \mbox{in}\
\mathcal{O},\vspace{2mm} \\
h\cdot \mathrm{n}=0,\qquad \left[ 2(\mathrm{n}\cdot Dh)+\gamma h\right]
\cdot \mathrm{\tau }=0 & \quad \mbox{on}\ \Gamma \vspace{2mm}.%
\end{array}%
\right.  \label{GS_NS}
\end{equation}
Next, we state a lemma concerning the regularity properties of the solution
of this system.

\begin{lemma}
\label{L_SP} Suppose $f\in H^{m}(\mathcal{O})$, $m=0,1$. Then system (\ref%
{GS_NS}) has a solution $(h,p)\in H^{m+2}(\mathcal{O})\times H^{m+1}(%
\mathcal{O})$, moreover the following estimates hold
\begin{align}
\Vert h\Vert _{H^{2}}& \leq C\Vert f\Vert _{2},  \label{SP_L2} \\
\Vert h\Vert _{H^{3}}& \leq C\Vert f\Vert _{H^{1}}.  \label{SP_L3}
\end{align}
\end{lemma}

\textbf{Proof.} Supposing that $f\in L^{2}(\mathcal{O})$, the existence of
the solution $(h,p)$ with $h$ in $H^{1}(\mathcal{O})$ is given by the
Lax-Millgram lemma. Multiplying (\ref{GS_NS})$_{1}$ by $h,$ \ we derive
\begin{equation*}
\left\Vert h\right\Vert _{2}^{2}+\alpha \left( 2\left\Vert Dh\right\Vert
_{2}^{2}+\gamma \Vert h\Vert _{L^{2}(\Gamma )}^{2}\right) =(f,h)\leq
\left\Vert f\right\Vert _{2}\left\Vert h\right\Vert _{2},
\end{equation*}%
which gives
\begin{equation}
\left\Vert h\right\Vert _{H^{1}}\leq C\left\Vert f\right\Vert _{2}.
\label{upo}
\end{equation}%
On the other hand, applying the operator $\mathrm{curl}$ to system (\ref%
{GS_NS}), we derive the following system for $u=\mathrm{curl}\,h$
\begin{equation}
\left\{
\begin{array}{ll}
u-\alpha \Delta u=\mathrm{curl}\,f\qquad & \quad \mbox{in}\ \mathcal{O},%
\vspace{2mm} \\
u=g(h)=\left( 2k-\gamma \right) \,h\cdot \mathrm{\tau } & \quad \mbox{on}\
\Gamma .%
\end{array}%
\right.  \label{ec}
\end{equation}%
Let us denote the extension of the unit exterior normal $\mathrm{n}$ (and
the tangent $\mathrm{\tau }=(-n_{2},n_{1})$) on the whole domain $\overline{%
\mathcal{O}}$ by the same notation $\mathrm{n}$ (and $\mathrm{\tau }$). \
Then the function $z=u-\left( 2k-\gamma \right) \,h\cdot \mathrm{\tau }$
solves the system
\begin{equation}
\left\{
\begin{array}{ll}
z-\alpha \Delta z=\mathrm{curl}\,f-\left( 2k-\gamma \right) \,h\cdot \mathrm{%
\tau }+\alpha \Delta \left[ \left( 2k-\gamma \right) \,h\cdot \mathrm{\tau }%
\right] \qquad & \quad \mbox{in}\ \mathcal{O},\vspace{2mm} \\
z=0 & \quad \mbox{on}\ \Gamma .%
\end{array}%
\right.  \label{ecz}
\end{equation}%
Multiplying equation (\ref{ecz})$_{1}$ by $z,$ integrating by parts and
using (\ref{upo}), we deduce
\begin{equation*}
\Vert z\Vert _{2}+\alpha \Vert \nabla z\Vert _{2}\leq C\left( \Vert f\Vert
_{2}+\Vert h\Vert _{H^{1}}\right) \leq C\Vert f\Vert _{2},
\end{equation*}%
which implies
\begin{equation}
\Vert u\Vert _{H^{1}}\leq C\left( \Vert f\Vert _{2}+\Vert h\Vert
_{H^{1}}\right) \leq C\Vert f\Vert _{2}.  \label{28}
\end{equation}%
In addition estimate $(2.3.3.7)$, p. 110 of \cite{G85} for system (\ref{ec})
gives
\begin{eqnarray}
\Vert u\Vert _{H^{2}} &\leq &C\left( \Vert \mathrm{curl}\,f\Vert
_{2}+\left\Vert \left( 2k-\gamma \right) \,h\cdot \mathrm{\tau }\right\Vert
_{H^{2-\frac{1}{2}}(\Gamma )}\right)  \notag \\
&\leq &C\left( \Vert f\Vert _{H^{1}}+\left\Vert h\right\Vert _{H^{2}}\right)
.  \label{jxx1}
\end{eqnarray}%
Since $h$ solves system (\ref{GS_NS}), then there exists a stream function $%
\varphi $ such that $h=\nabla ^{\bot }\varphi ,$ satisfying the system
\begin{equation}
\left\{
\begin{array}{ll}
\Delta \varphi =u\qquad & \quad \mbox{in}\ \mathcal{O},\vspace{2mm} \\
\varphi =0 & \quad \mbox{on}\ \Gamma%
\end{array}%
\right.  \label{es}
\end{equation}%
and the estimate
\begin{equation}
\Vert \varphi \Vert _{H^{2+m}}\leq \Vert u\Vert _{H^{m}},\qquad m\in \mathbb{%
N}_{0},  \label{esf}
\end{equation}%
by Theorem 2.5.1.1 , p. 128 of \cite{G85}.

Combining (\ref{28}) and (\ref{esf}) with $m=1$, we deduce
\begin{equation*}
\Vert \varphi \Vert _{H^{3}}\leq C\Vert u\Vert _{H^{1}}\leq C\Vert f\Vert
_{2},
\end{equation*}%
hence $h=\nabla ^{\bot }\varphi \in H^{2}$ and (\ref{SP_L2}) hold. \
Moreover (\ref{jxx1}) and (\ref{esf}) with $m=2$ imply
\begin{equation*}
\Vert \varphi \Vert _{H^{4}}\leq \Vert u\Vert _{H^{2}}\leq C\left( \Vert
f\Vert _{H^{1}}+\Vert h\Vert _{H^{2}}\right) .
\end{equation*}%
Invoking (\ref{SP_L2}), we conclude that $h=\nabla ^{\bot }\varphi \in
H^{3} $ and (\ref{SP_L3}) hold.\bigskip $\hfill \hfill \blacksquare \hfill $

\bigskip

Let us recall that the space $W$ introduced in (\ref{w}) is naturally
endowed with the Sobolev norm $\Vert \cdot \Vert _{H^{2}}$. The next result
follows directly from Lemma 5 in \cite{BI06} and helps to introduce on $W$
an equivalent norm that will be useful to analyze the stability in Section 5.

\begin{lemma}
\label{sigma_psigma} For each $y\in W$, we have
\begin{equation}
\left\Vert \upsilon (y)-\mathbb{P}\upsilon (y)\right\Vert _{2}\leq
C\left\Vert y\right\Vert _{H^{1}},  \label{sigma_psigma1}
\end{equation}%
\begin{equation}
\left\Vert \upsilon (y)-\mathbb{P}\upsilon (y)\right\Vert _{H^{1}}\leq
C\left\Vert y\right\Vert _{H^{2}}.  \label{sigma_psigma2}
\end{equation}
\end{lemma}

The next regularity result will be fundamental to establish the
well-posedness of the velocity equation (see Propositions 6 in \cite{BR03}
and Lemma 2.1 in \cite{CG97} for similar results).

\begin{lemma}
\label{non_lin}Let $y\in \widetilde{W}$. Then, the following estimates hold%
\begin{align}
\left\Vert y\right\Vert _{H^{2}}& \leq C\left( \left\Vert \mathbb{P}%
\,\upsilon (y)\right\Vert _{2}+\left\Vert y\right\Vert _{H^{1}}\right) ,
\label{www} \\
\left\Vert y\right\Vert _{H^{3}}& \leq C\left( \left\Vert \mathrm{curl}%
\,\upsilon (y)\right\Vert _{2}+\left\Vert y\right\Vert _{H^{1}}\right) .
\label{y2_sigma}
\end{align}
\end{lemma}

\textbf{Proof.} Considering system (\ref{GS_NS}) with $f=\upsilon (y),$ then
the pair $(y,0)$ is obviously the solution of such system. Hence estimate (%
\ref{SP_L2}) yields%
\begin{equation*}
\left\Vert y\right\Vert _{H^{2}}\leq C\left\Vert \upsilon (y)\right\Vert
_{2}\leq C\left( \left\Vert \upsilon (y)-\mathbb{P}\upsilon (y)\right\Vert
_{2}+\left\Vert \mathbb{P}\upsilon (y)\right\Vert _{2}\right) .
\end{equation*}%
Applying (\ref{sigma_psigma1}), we deduce (\ref{www}).

Since $~\mathrm{curl}\,\upsilon (y)\in L^{2}(\mathcal{O})~$ and $~\nabla
\cdot \left( \mathrm{curl}\,\upsilon (y)\right) =0$,$~$ there exists a
unique vector-potential $\psi \in H^{1}(\mathcal{O})$ such that
\begin{equation*}
\left\{
\begin{array}{ll}
\mathrm{{curl}\,}\psi =\mathrm{curl}\,\upsilon (y),\qquad \mathrm{div}\,\psi
=0 & \quad \mbox{in}\ \mathcal{O},\vspace{2mm} \\
\psi \cdot \mathrm{n}=0 & \quad \mbox{on}\ \Gamma%
\end{array}%
\right.
\end{equation*}%
and
\begin{equation}
\left\Vert \psi \right\Vert _{H^{1}}\leq C\left\Vert \mathrm{curl}\,\upsilon
(y)\right\Vert _{2}.  \label{sigma_phi}
\end{equation}%
It follows that $\mathrm{curl}\left( y-\alpha \Delta y-\psi \right) =0 $ and
there exists $~\pi \in L^{2}(\mathcal{O}),~$ such that~%
\begin{equation*}
y-\alpha \Delta y-\psi +\nabla \pi =0.
\end{equation*}%
Hence $y$ is the solution of the Stokes system (\ref{GS_NS}) where $f$ is
replaced by $\psi .$

As a consequence of (\ref{SP_L3}), we have
\begin{equation}
\left\Vert y\right\Vert _{H^{3}}\leq C\left( \left\Vert \psi \right\Vert
_{H^{1}}+\left\Vert y\right\Vert _{H^{1}}\right) .  \label{y_phi}
\end{equation}%
Using (\ref{sigma_phi}) we obtain the claimed result (\ref{y2_sigma}%
).\bigskip $\hfill \hfill \blacksquare \hfill $

\bigskip

In order to define the solution of \ equation (\ref{equation_etat})$_{1}$ in
the distributional sense, we introduce a trilinear functional that is well
known in the context of the Navier-Stokes equations
\begin{equation*}
b(\phi ,z,y)=\left( \phi \cdot \nabla z,y\right) ,\qquad \forall \phi
,z,y\in V.
\end{equation*}%
In what follows we often will use the following property
\begin{equation}
b(\phi ,z,y)=-b(\phi ,y,z),  \label{ra}
\end{equation}%
that follows taking integration by parts, knowing that $\phi $ is divergence
free and $(\phi \cdot n)=0$ on $\Gamma $.

Straightforward computations yield the following relation \label%
{non_lin_curl}
\begin{equation}
\left( \mathrm{curl}\,\upsilon (y)\times z,\phi \right) =b\left( \phi
,z,\upsilon (y)\right) -b\left( z,\phi ,\upsilon (y)\right) \quad \quad
\forall y\in \widetilde{W},\;z,\phi \in V.  \label{bb1}
\end{equation}

In the next lemma, we deduce crucial estimates of major importance to
establish \ the well-posedness of system (\ref{equation_etat}), as well as
to prove the stability property of their solutions. We should mention that
some estimates follow from an adaptation of the method considered in \cite%
{BR03} to prove the uniqueness.

\begin{lemma}
\label{rm2} Let $y,$ $z,\;\phi \in \widetilde{W}$. \ Then
\begin{align}
\left\vert \left( \mathrm{curl}\,\upsilon (y)\times z,\phi \right)
\right\vert & \leq C\left\Vert y\right\Vert _{H^{3}}\left\Vert z\right\Vert
_{H^{1}}\left\Vert \phi \right\Vert _{H^{3}},  \label{rm2_lin1} \\
\left\vert \left( \mathrm{curl}\,\upsilon (y)\times z,\phi \right)
\right\vert & \leq C\left\Vert y\right\Vert _{H^{1}}\left\Vert z\right\Vert
_{H^{3}}\left\Vert \phi \right\Vert _{H^{3}},  \label{rm2_lin3} \\
\left\vert \left( \mathrm{curl}\,\upsilon (y)\times z,y\right) \right\vert &
\leq C\left\Vert y\right\Vert _{H^{1}}^{2}\left\Vert z\right\Vert _{H^{3}}.
\label{rm2_lin}
\end{align}
\end{lemma}

\bigskip \textbf{Proof. }\textit{1st step. The proof of estimate (\ref%
{rm2_lin1}). }We directly can estimate
\begin{equation*}
\left\vert \left( \mathrm{curl}\,\upsilon (y)\times z,\phi \right)
\right\vert \leq \Vert \phi \Vert _{\infty }\Vert z\Vert _{2}\Vert \mathrm{%
curl}\,\upsilon (y)\Vert _{2}\leq \Vert \phi \Vert _{H^{3}}\Vert z\Vert
_{2}\Vert y\Vert _{H^{3}}
\end{equation*}%
by Sobolev's embedding $H^{3}(\mathcal{O})\hookrightarrow L^{\infty }(%
\mathcal{O}).$ Hence we have (\ref{rm2_lin1}).

\textit{2nd step.} \textit{The proof of estimate (\ref{rm2_lin3}).\ \ \ }%
Equality (\ref{bb1}) gives%
\begin{equation}
\left( \mathrm{curl}\,\upsilon (y)\times z,\phi \right) =b(\phi
,z,y)-b(z,\phi ,y)-\alpha \left( b(\phi ,z,\Delta y)-b(z,\phi ,\Delta
y)\right) .  \label{mk}
\end{equation}%
With the help of Sobolev's embedding $H^{1}(\mathcal{O})\hookrightarrow
L^{4}(\mathcal{O})$, it is easy to see that
\begin{align}
\left\vert b\left( \phi ,z,y\right) -b\left( z,\phi ,y\right) \right\vert &
=\left\vert b\left( \phi ,z,y\right) +b\left( z,\phi ,y\right) \right\vert
\notag \\
& \leq \Vert \phi \Vert _{4}\Vert \nabla z\Vert _{2}\Vert y\Vert _{4}+\Vert
z\Vert _{4}\Vert \nabla \phi \Vert _{2}\Vert y\Vert _{4}  \notag \\
& \leq C\Vert \phi \Vert _{H^{1}}\Vert z\Vert _{H^{1}}\Vert y\Vert _{H^{1}}.
\label{10-2-1}
\end{align}%
Integrating by parts and using the boundary conditions, we derive
\begin{align}
b\left( \phi ,z,\Delta y\right) & =\sum_{i,j=1}^{2}\int_{\mathcal{O}}\phi
_{i}\tfrac{\partial z_{j}}{\partial x_{i}}\Delta
y_{j}=\sum_{i,j,k=1}^{2}\int_{\mathcal{O}}\phi _{i}\tfrac{\partial z_{j}}{%
\partial x_{i}}\tfrac{\partial }{\partial x_{k}}\left( \tfrac{\partial y_{j}%
}{\partial x_{k}}-\tfrac{\partial y_{k}}{\partial x_{j}}\right)
=\sum_{i,j,k=1}^{2}\int_{\mathcal{O}}\phi _{i}\tfrac{\partial z_{j}}{%
\partial x_{i}}\tfrac{\partial }{\partial x_{k}}(A_{jk}(y))  \notag \\
& =\sum_{i,j,k=1}^{2}\int_{\Gamma }\phi _{i}\tfrac{\partial z_{j}}{\partial
x_{i}}A_{jk}(y)n_{k}-\sum_{i,j,k=1}^{2}\int_{\mathcal{O}}\tfrac{\partial }{%
\partial x_{k}}\left( \phi _{i}\tfrac{\partial z_{j}}{\partial x_{i}}\right)
A_{jk}(y)  \notag \\
& =\sum_{i,j=1}^{2}\int_{\Gamma }\phi _{i}\tfrac{\partial z_{j}}{\partial
x_{i}}\,g(y)\tau _{j}-\sum_{i,j,k=1}^{2}\int_{\mathcal{O}}\tfrac{\partial
\phi _{i}}{\partial x_{k}}\tfrac{\partial z_{j}}{\partial x_{i}}%
A_{jk}(y)-\sum_{i,j,k=1}^{2}\int_{\mathcal{O}}\phi _{i}\tfrac{\partial
^{2}z_{j}}{\partial x_{k}\partial x_{i}}A_{jk}(y)  \notag \\
& =I_{1}+I_{2}+I_{3}.  \label{b}
\end{align}%
Again, integrating by parts, it follows that
\begin{align*}
I_{1}& =\sum_{i,j=1}^{2}\int_{\Gamma }\phi _{i}\tfrac{\partial z_{j}}{%
\partial x_{i}}\,g(y)\tau _{j}=\sum_{i,j=1}^{2}\int_{\mathcal{O}}\phi _{i}%
\tfrac{\partial z_{j}}{\partial x_{i}}\,g(y)\tau _{j}\,\mathrm{div}\,\mathrm{%
n}+\sum_{i,j,k=1}^{2}\int_{\mathcal{O}}\tfrac{\partial }{\partial x_{k}}%
\bigl(\phi _{i}\tfrac{\partial z_{j}}{\partial x_{i}}\,g(y)\tau _{j}\bigr)%
n_{k} \\
& =b\left( \phi ,z,\mathrm{div}\,\mathrm{n}\,g(y)\mathrm{\tau }\right)
+b\left( \left( \mathrm{n}\cdot \nabla \right) \phi ,z,\,n_{k}\,g(y)\mathrm{%
\tau }\right) +\sum_{k=1}^{2}b\left( \phi ,\,\tfrac{\partial z}{\partial
x_{k}},n_{k}\,g(y)\mathrm{\tau }\right) \\
& +b\left( \phi ,z,\,\left( \mathrm{n}\cdot \nabla \right) \bigl(g(y)\mathrm{%
\tau }\bigr)\right) .
\end{align*}%
Then, using Sobolev's embedding $H^{2}(\mathcal{O})\hookrightarrow L^{\infty
}(\mathcal{O}),$ we easily derive
\begin{align}
& \left\vert b\left( \phi ,z,\Delta y\right) \right\vert \leq \left\vert
b\left( \phi ,z,\mathrm{div}\,\mathrm{n}\,g(y)\mathrm{\tau }\right)
\right\vert +\left\vert b\left( \left( \mathrm{n}\cdot \nabla \right) \phi
,z,\,\,g(y)\mathrm{\tau }\right) \right\vert +\left\vert b\left( \phi
,z,\,\left( \mathrm{n}\cdot \nabla \right) \bigl(g(y)\mathrm{\tau }\bigr)%
\right) \right\vert  \notag \\
& +\sum_{k=1}^{2}\left( \left\vert b\left( \phi ,\tfrac{\partial z}{\partial
x_{k}},\,n_{k}\,g(y)\mathrm{\tau }\right) \right\vert +\left\vert b\left(
\tfrac{\partial \phi }{\partial x_{k}},z,A_{\cdot \,k}(y)\right) \right\vert
+\left\vert b\left( \phi ,\tfrac{\partial z}{\partial x_{k}},A_{\cdot
\,k}(y)\right) \right\vert \right)  \notag \\
& \leq C\Vert z\Vert _{H^{3}}\Vert \phi \Vert _{H^{3}}\Vert y\Vert _{H^{1}}.
\label{trilin_100}
\end{align}%
By symmetry, it follows that
\begin{equation}
\left\vert b\left( z,\phi ,\Delta y\right) \right\vert \leq C\Vert z\Vert
_{H^{3}}\Vert \phi \Vert _{H^{3}}\Vert y\Vert _{H^{1}}.  \label{lk}
\end{equation}%
Then (\ref{rm2_lin3}) follows from (\ref{mk})-(\ref{lk}).

\textit{3d step.} \textit{The proof of estimate (\ref{rm2_lin}). \ } As in
the above computations for (\ref{b}), we obtain
\begin{align*}
b\left( z,y,\Delta y\right) & =\sum_{i,j=1}^{2}\int_{\Gamma }z_{i}\tfrac{%
\partial y_{j}}{\partial x_{i}}\,g(y)\tau _{j}-\sum_{i,j,k=1}^{2}\int_{%
\mathcal{O}}\tfrac{\partial z_{i}}{\partial x_{k}}\tfrac{\partial y_{j}}{%
\partial x_{i}}A_{jk}(y)-\sum_{i,j,k=1}^{2}\int_{\mathcal{O}}z_{i}\tfrac{%
\partial ^{2}y_{j}}{\partial x_{k}\partial x_{i}}A_{jk}(y) \\
& =J_{1}+J_{2}+J_{3},
\end{align*}%
where
\begin{align*}
J_{1}& =b\left( z,y,\mathrm{div}\,\mathrm{n}\,g(y)\mathrm{\tau }\right)
+b\left( \left( \mathrm{n}\cdot \nabla \right) z,y,\,g(y)\mathrm{\tau }%
\right) +\sum_{k=1}^{2}b\left( z,\,\tfrac{\partial y}{\partial x_{k}}%
,n_{k}\,g(y)\mathrm{\tau }\right) \\
& +b\left( z,y,\left( \mathrm{n}\cdot \nabla \right) \bigl(g(y)\mathrm{\tau }%
\bigr)\right) , \\
J_{2}& =-\sum_{k=1}^{2}b\left( \tfrac{\partial z}{\partial x_{k}},y,A_{\cdot
\,k}(y)\right)
\end{align*}%
and%
\begin{eqnarray*}
J_{3} &=&-\sum_{i,j,k=1}^{2}\int_{\mathcal{O}}z_{i}\tfrac{\partial }{%
\partial x_{i}}\left( \tfrac{\partial y_{j}}{\partial x_{k}}\right) \left(
\tfrac{\partial y_{j}}{\partial x_{k}}-\tfrac{\partial y_{k}}{\partial x_{j}}%
\right) =\sum_{i,j,k=1}^{2}\int_{\mathcal{O}}z_{i}\tfrac{\partial }{\partial
x_{i}}\left( \tfrac{\partial y_{j}}{\partial x_{k}}\right) \tfrac{\partial
y_{k}}{\partial x_{j}} \\
&=&\frac{1}{2}\sum_{i,j,k=1}^{2}\int_{\mathcal{O}}z_{i}\left[ \tfrac{%
\partial }{\partial x_{i}}\left( \tfrac{\partial y_{j}}{\partial x_{k}}%
\right) \tfrac{\partial y_{k}}{\partial x_{j}}+\tfrac{\partial }{\partial
x_{i}}\left( \tfrac{\partial y_{k}}{\partial x_{j}}\right) \tfrac{\partial
y_{j}}{\partial x_{k}}\right] =\sum_{i,j,k=1}^{2}\int_{\mathcal{O}}z_{i}%
\tfrac{\partial }{\partial x_{i}}\left( \tfrac{\partial y_{j}}{\partial x_{k}%
}\tfrac{\partial y_{k}}{\partial x_{j}}\right) =0.
\end{eqnarray*}%
Therefore, we derive
\begin{align}
\left\vert b\left( z,y,\Delta y\right) \right\vert & \leq \left\vert b\left(
z,y,\mathrm{div}\,\mathrm{n}\,g(y)\mathrm{\tau }\right) \right\vert
+\left\vert b\left( \left( \mathrm{n}\cdot \nabla \right) z,y,\,g(y)\mathrm{%
\tau }\right) \right\vert +\sum_{k=1}^{2}\left\vert b\left( z,\,n_{k}\,g(y)%
\mathrm{\tau },\tfrac{\partial y}{\partial x_{k}}\right) \right\vert  \notag
\\
& +\left\vert b\left( z,y,\left( \mathrm{n}\cdot \nabla \right) \bigl(g(y)%
\mathrm{\tau }\bigr)\right) \right\vert +\sum_{k=1}^{2}\left\vert b\left(
\tfrac{\partial z}{\partial x_{k}},y,A_{\cdot \,k}(y)\right) \right\vert
\leq C\Vert z\Vert _{H^{2}}\Vert y\Vert _{H^{1}}^{2},  \label{trilin_1}
\end{align}%
where we have used that
\begin{equation*}
\sum_{k=1}^{2}b\left( z,\,\tfrac{\partial y}{\partial x_{k}},n_{k}\,g(y)%
\mathrm{\tau }\right) =-\sum_{k=1}^{2}b\left( z,\,n_{k}\,g(y)\mathrm{\tau },%
\tfrac{\partial y}{\partial x_{k}}\right)
\end{equation*}%
by (\ref{ra}).

Taking $\ \phi =y$ \ in (\ref{b}) we have
\begin{equation*}
b\left( y,z,\Delta y\right) =\sum_{i,j=1}^{2}\int_{\Gamma }y_{i}\tfrac{%
\partial z_{j}}{\partial x_{i}}g(y)\tau _{j}-\sum_{k=1}^{2}b\left( \tfrac{%
\partial y}{\partial x_{k}},z,A_{\cdot \,k}(y)\right) -\sum_{k=1}^{2}b\left(
y,\tfrac{\partial z}{\partial x_{k}},A_{\cdot \,k}(y)\right) .
\end{equation*}%
%
%
%
%
%
%
%
%
%
%
%
%
%
%
%
%
%
%
%
%
%
%
%
%
%
%
%
%
%
%
%
%
%
%
%
%
%
%
%
%
%
Taking into account the embedding theorems $H^{2}(\mathcal{O}%
)\hookrightarrow C(\overline{\mathcal{O}}),$ $H^{1}(\mathcal{O}%
)\hookrightarrow L^{2}(\Gamma )$ \ and $H^{1}(\mathcal{O})\hookrightarrow
L^{4}(\mathcal{O}),$ we have%
\begin{align}
|b\left( y,z,\Delta y\right) |& \leq \sum_{i,j=1}^{2}\int_{\Gamma }\bigl|%
y_{i}\tfrac{\partial z_{j}}{\partial x_{i}}g(y)\tau _{j}\bigr|%
+\sum_{k=1}^{2}\left( \bigl|b\left( \tfrac{\partial y}{\partial x_{k}}%
,z,A_{\cdot \,k}(y)\right) \bigr|+\bigl|b\left( y,\tfrac{\partial z}{%
\partial x_{k}},A_{\cdot \,k}(y)\right) \bigr|\right)  \notag \\
& \leq C\Vert y\Vert _{L^{2}(\Gamma )}^{2}\left\Vert \nabla z\right\Vert _{C(%
\overline{\mathcal{O}})}+C\left\Vert \nabla y\right\Vert _{2}^{2}\left\Vert
\nabla z\right\Vert _{\infty }+\sum_{i,j=1}^{2}\left\Vert y\right\Vert
_{4}\left\Vert \tfrac{\partial ^{2}z}{\partial x_{i}\partial x_{k}}%
\right\Vert _{4}\left\Vert \nabla y\right\Vert _{2}\vspace{2mm}  \notag \\
& \leq C\left\Vert y\right\Vert _{H^{1}}^{2}\left\Vert z\right\Vert _{H^{3}}.
\label{GF_2mm}
\end{align}%
Then (\ref{rm2_lin}) is a consequence of (\ref{10-2-1}) and (\ref{trilin_1}%
)-(\ref{GF_2mm}).$\hfill \hfill \blacksquare \hfill $

\bigskip


\section{Existence of strong solution}

\setcounter{equation}{0}

\bigskip

The aim of the present section is to establish the existence of a strong
solution for system (\ref{equation_etat}) in the probabilistic sense.

\begin{definition}
\label{4.1} Let
\begin{equation*}
U\in L^{2}(\Omega \times (0,T);H(\mathrm{curl};\mathcal{O})),\text{\quad
\quad }Y_{0}\in L^{2}(\Omega ,\widetilde{W}).
\end{equation*}%
A stochastic process $Y\in L^{2}(\Omega ,L^{\infty }(0,T;\widetilde{W}))$ is
a strong solution of $(\ref{equation_etat})$, if \ for a.e.-$P$ and \ a.e. $%
t\in (0,T),$ the following equation holds
\begin{align}
\left( \upsilon (Y(t)),\phi \right) =& \int_{0}^{t}\left[ -2\nu \left(
DY(s),D\phi \right) -\nu \gamma \int_{\Gamma }y\cdot \phi \ dx-\left(
\mathrm{curl}\,\upsilon (Y(s))\times Y(s),\phi \right) \right] \,ds  \notag
\\
& +\left( \upsilon (Y(0)),\phi \right) +\int_{0}^{t}\left( U(s),\phi \right)
\,ds+\int_{0}^{t}\left( G(s,Y(s)),\phi \right) \,dW_{s}
\label{var_form_state}
\end{align}%
for all $\phi \in V$, where the nonlinear term should be understood in the
sense
\begin{equation*}
\begin{array}{ll}
\left( \mathrm{curl}\,\upsilon (Y(t))\times Y(t),\phi \right) & =b\left(
\phi ,Y(t),\upsilon (Y(s))\right) -b\left( Y(t),\phi ,\upsilon (Y(s))\right)
\vspace{2mm}
\end{array}%
\end{equation*}
and the stochastic integral is defined by
$$
\int_{0}^{t}\left( G(s,Y(s)),\phi \right) \,dW_{s}=\sum_{k=1}^m \int_{0}^{t} \left( G^k(s,Y(s)),\phi \right) \,dW^k_{s}.
$$
\vspace{2mm}\newline
\end{definition}

Let us formulate our main existence and uniqueness result, which will be
shown in this section.

\begin{theorem}
\label{the_1} Assume that
\begin{equation*}
U\in L^{p}(\Omega \times (0,T);H(\mathrm{curl};\mathcal{O})),\qquad Y_{0}\in
L^{p}(\Omega ,V)\cap L^{2}(\Omega ,\widetilde{W})\qquad \text{for some\ }\
4\leq p<\infty .
\end{equation*}%
Then there exists a unique solution $Y$ to equation $(\ref{var_form_state})$
which belongs to
\begin{equation*}
L^{2}(\Omega ,L^{\infty }(0,T;\widetilde{W}))\cap L^{p}(\Omega ,L^{\infty
}(0,T;V)).
\end{equation*}%
Moreover, the following estimates hold
\begin{equation*}
\frac{1}{2}\mathbb{E}\sup_{s\in \lbrack 0,t]}\left\Vert Y(s)\right\Vert
_{V}^{2}+\mathbb{E}\int_{0}^{t}\left( 4\nu \left\Vert DY\right\Vert
_{2}^{2}\,+2\nu \gamma \left\Vert Y\right\Vert _{L^{2}(\Gamma )}^{2}\right)
\,ds\leq C\left( \mathbb{E}\left\Vert Y_{0}\right\Vert _{V}^{2}+\mathbb{E}%
\Vert U\Vert _{L^{2}(0,t;L^{2})}^{2}+1\right) ,
\end{equation*}%
\begin{equation*}
\mathbb{E}\sup_{s\in \lbrack 0,t]}\left\Vert \mathrm{curl}\,\upsilon \left(
Y(s)\right) \right\Vert _{2}^{2}\leq C\left( \mathbb{E}\left\Vert \mathrm{%
curl}\,\upsilon (Y_{0})\right\Vert _{2}^{2}+\mathbb{E}\left\Vert
U\right\Vert _{L^{2}(0,t;H^{1})}^{2}+1\right) .
\end{equation*}
\end{theorem}

\bigskip

The proof of the theorem is given by \textit{Galerkin's approximation method}%
. We consider the inner product of $\widetilde{W}$ defined by
\begin{equation}
\left( y,z\right) _{\widetilde{W}}=\left( \mathrm{curl}\upsilon (y),\mathrm{%
curl}\upsilon (z)\right) +\left( y,z\right) _{V}.  \label{inpr}
\end{equation}%
Taking into account (\ref{korn}) and (\ref{y2_sigma}) the norm $\Vert \cdot
\Vert _{\widetilde{W}}$ induced by this inner product is equivalent to $%
\Vert \cdot \Vert _{H^{3}}$. The injection operator $I:\widetilde{W}%
\rightarrow V$ is a compact operator, then there exists a basis $%
\{e_{i}\}\subset \widetilde{W}$ of eigenfunctions
\begin{equation}
\left( y,e_{i}\right) _{\widetilde{W}}=\lambda _{i}\left( y,e_{i}\right)
_{V},\qquad \forall y\in \widetilde{W},\;i\in \mathbb{N},  \label{y4}
\end{equation}%
being  an orthonormal basis for $V$ and the corresponding
sequence $\{\lambda _{i}\}$ of eigenvalues verifies $\lambda _{i}>0$, $%
\forall i\in \mathbb{N}$ and $\lambda _{i}\rightarrow \infty $ as $%
i\rightarrow \infty .$ Let us notice that the ellipticity of equation (\ref%
{y4}) increases the regularity of their solutions. Hence without loss of
generality we can consider $\{e_{i}\}\subset H^{4}$ (see \cite{BR03}).
\vspace{2mm}\newline
In this section, we consider this basis and introduce the Faedo-Galerkin
approximation of system (\ref{equation_etat}). Let $W_{n}=\mathrm{span}%
\,\{e_{1},\ldots ,e_{n}\}$ and define
\begin{equation*}
Y_{n}(t)=\sum_{j=1}^{n}c_{j}^{n}(t)e_{j}
\end{equation*}%
as the solution of the stochastic differential equation
\begin{equation}
\left\{
\begin{array}{l}
d\left( \upsilon \left( Y_{n}\right) ,\phi \right) =\left( \left( \nu \Delta
Y_{n}-\mathrm{curl}\left( \upsilon \left( Y_{n}\right) \right) \times
Y_{n}+U\right) ,\phi \right) \,dt+\left( G(t,Y_{n}),\phi \right) \,dW_{t},%
\vspace{2mm} \\
Y_{n}(0)=Y_{n,0},\qquad \qquad \qquad \qquad \qquad \qquad \forall \phi \in
W_{n}.%
\end{array}%
\right.  \label{y1}
\end{equation}%
Here $Y_{n,0}$ denotes the projection of the initial condition $Y_{0}$ onto
the space $W_{n}$.

Let us notice that $\{\widetilde{e}_{j}=\frac{1}{\sqrt{\lambda _{j}}}%
e_{j}\}_{j=1}^{\infty }$ is an orthonormal basis for $\widetilde{W}$ and
\begin{equation*}
Y_{n,0}=\sum_{j=1}^{n}\left( Y_{0},e_{j}\right)
_{V}e_{j}=\sum_{j=1}^{n}\left( Y_{0},\widetilde{e}_{j}\right) _{\widetilde{W}%
}\widetilde{e}_{j},
\end{equation*}%
then the Parseval's identity gives
\begin{equation}
\left\Vert Y_{n}(0)\right\Vert _{V}\leq \left\Vert Y_{0}\right\Vert
_{V}\quad \text{ and }\quad \left\Vert Y_{n}(0)\right\Vert _{\widetilde{W}%
}\leq \left\Vert Y_{0}\right\Vert _{\widetilde{W}}.  \label{CI}
\end{equation}%
Equation (\ref{y1}) defines a system of stochastic ordinary differential
equations in $\mathbb{R}^{n}$ with locally Lipschitz nonlinearities. Hence
there exists a local-in-time solution $\ Y_{n}$ as an adapted process in the
space $C([0,T_{n}];W_{n})$. \ The global-in-time existence of $Y_{n}$
follows from uniform estimates on $n=1,2,...,$ that will be deduced in the
next Lemma \ (a similar reasoning can be found in \cite{ABW10}, \cite{RS12}%
). 

\begin{lemma}
\label{existence_state} Assume that%
\begin{equation*}
U\in L^{2}(\Omega \times (0,T);H(\mathrm{curl};\mathcal{O})),\qquad Y_{0}\in
L^{2}(\Omega ,\widetilde{W}).
\end{equation*}%
Then problem (\ref{y1}) admits a unique solution $Y_{n}\in L^{2}(\Omega
,L^{\infty }(0,T;\widetilde{W}))$. \ Furthermore, for any $t\in \lbrack 0,T]$%
, the following estimates hold
\begin{eqnarray}
\frac{1}{2}\mathbb{E}\sup_{s\in \lbrack 0,t]}\left\Vert Y_{n}(s)\right\Vert
_{V}^{2} &+&\mathbb{E}\int_{0}^{t}\left( 4\nu \left\Vert DY_{n}\right\Vert
_{2}^{2}\,+2\nu \gamma \left\Vert Y_{n}\right\Vert _{L^{2}(\Gamma
)}^{2}\right) \,ds  \notag \\
&\leq &C\left( 1+\mathbb{E}\left\Vert Y_{0}\right\Vert _{V}^{2}+\mathbb{E}%
\Vert U\Vert _{L^{2}(0,t;L^{2})}^{2}\right) ,  \label{ineq1}
\end{eqnarray}%
\begin{align}
\frac{1}{2}\mathbb{E}\sup_{s\in \lbrack 0,t]}\left\Vert \mathrm{curl}%
\,\upsilon \left( Y_{n}(s)\right) \right\Vert _{2}^{2}& +\frac{2\nu }{\alpha
}\mathbb{E}\int_{0}^{t}\left\Vert \mathrm{curl}\,\upsilon \left(
Y_{n}\right) \right\Vert _{2}^{2}\,ds\leq \mathbb{E}\left\Vert \mathrm{curl}%
\,\upsilon \left( Y_{0}\right) \right\Vert _{2}^{2}  \notag \\
& +C\mathbb{E}\int_{0}^{t}\left\Vert \mathrm{curl}\,U\right\Vert
_{2}^{2}\,ds+C\mathbb{E}\int_{0}^{t}\left( 1+\left\Vert Y_{n}\right\Vert
_{V}^{2}\right) \,ds  \label{ineq222}
\end{align}%
and%
\begin{equation}
\mathbb{E}\sup_{s\in \lbrack 0,t]}\left\Vert Y_{n}(s)\right\Vert _{%
\widetilde{W}}^{2}\leq C(\mathbb{E}\left\Vert Y_{0}\right\Vert _{\widetilde{W%
}}^{2}+\mathbb{E}\left\Vert \,U\right\Vert _{L^{2}(0,t;H(\mathrm{curl};%
\mathcal{O}))}^{2}),  \label{ineq2}
\end{equation}%
where $C$ are positive constants independent of $n$ (and may depend on the
data of our problem the domain $\mathcal{O},$ \ the regularity of $\Gamma ,$
the physical constants $\nu $, $\alpha $, $\gamma ).$
\end{lemma}

\textbf{Proof.} {\ For each $n\in \mathbb{N}$, let us consider the sequence $%
\{\tau _{N}^{n}\}$}$_{{N\in \mathbb{N}}}${\ of the stopping times
\begin{equation*}
\tau _{N}^{n}=\inf \{t\geq 0:\Vert Y_{n}(t)\Vert _{H^{3}}\geq N\}\wedge
T_{n}.
\end{equation*}%
}

In order to simplify the notation, let us introduce the function
\begin{equation*}
f(Y_{n})=\left( \nu \Delta Y_{n}-\mathrm{curl}\left( \upsilon \left(
Y_{n}\right) \right) \times Y_{n}+U\right) \in H^{1}(\mathcal{O}).
\end{equation*}%
Taking $\phi =e_{i}$ for each $i=1,\dots ,n$\ in equation (\ref{y1}), we
obtain 

\begin{equation}
d\left( Y_{n},e_{i}\right) _{V}=\left( f(Y_{n}),e_{i}\right) \,dt+\left(
G(t,Y_{n}),e_{i}\right) \,dW_{t}.  \label{y3}
\end{equation}%
\textit{Step 1. Estimate in the space }$V$\textit{\ for }$Y_{n}$\textit{,
depending on the stopping times }$\tau _{N}^{n}.$

The It\^{o} formula gives
\begin{equation*}
d\left( Y_{n},e_{i}\right) _{V}^{2}=2\left( Y_{n},e_{i}\right) _{V}\left(
f(Y_{n}),e_{i}\right) \,dt+2\left( Y_{n},e_{i}\right) _{V}\left(
G(t,Y_{n}),e_{i}\right) \,dW_{t}+|\left( G\left( t,Y_{n}\right) ,e_{i}\right)|
^{2}\,dt,
\end{equation*}%
where the module in the last term is defined by \eqref{product}. Summing these equalities 
over $i=1,\dots ,n,$ we obtain 
\begin{equation*}
d\left\Vert Y_{n}\right\Vert _{V}^{2}=2\left( f(Y_{n}),Y_{n}\right)
\,dt+2\left( G(t,Y_{n}),Y_{n}\right) \,dW_{t}+\sum_{i=1}^{n}|\left( G\left(
t,Y_{n}\right) ,e_{i}\right)|^{2}\,dt.
\end{equation*}%
We know that
\begin{align}
\left( f(Y_{n}),Y_{n}\right) & =-2\nu \left\Vert DY_{n}\right\Vert
_{2}^{2}-\nu \gamma \left\Vert Y_{n}\right\Vert _{L^{2}(\Gamma )}^{2}+\alpha
\left( b\left( Y_{n},Y_{n},\Delta Y_{n}\right) -b\left( Y_{n},Y_{n},\Delta
Y_{n}\right) \right)  \notag \\
& -b\left( Y_{n},Y_{n},Y_{n}\right) +\left( U,Y_{n}\right)  \notag \\
& =-2\nu \left\Vert DY_{n}\right\Vert _{2}^{2}-\nu \gamma \left\Vert
Y_{n}\right\Vert _{L^{2}(\Gamma )}^{2}+\left( U,Y_{n}\right) ,
\label{equation_etat_temps_n1}
\end{align}%
hence%
\begin{eqnarray}
d\left\Vert Y_{n}\right\Vert _{V}^{2} &=&2\left( -2\nu \left\Vert
DY_{n}\right\Vert _{2}^{2}-\nu \gamma \left\Vert Y_{n}\right\Vert
_{L^{2}(\Gamma )}^{2}+\left( U,Y_{n}\right) \right) \,dt  \notag \\
&&+2\left( G(t,Y_{n}),Y_{n}\right) \,dW_{t}+\sum_{i=1}^{n} |\left( G\left(
t,Y_{n}\right) ,e_{i}\right)|^{2}\,dt, \label{sh1}
\end{eqnarray}

Let $\tilde{G}_{n}$ be the solution of (\ref{GS_NS}) for $%
f=G(t,Y_{n})$. Then
\begin{equation*}
(\tilde{G}_{n},e_{i})_{V}=(G(t,Y_{n}),e_{i})\quad \text{for }i=1,\dots ,n\
\end{equation*}%
which implies
\begin{equation}
\sum_{i=1}^{n}|\left( G\left( t,Y_{n}\right) ,e_{i}\right)| ^{2}=||\widetilde{G%
}_{n}||_{V}^{2}\leq C||G\left( t,Y_{n}\right) ||_{2}^{2}\leq
C(1+||Y_{n}||_{V}^{2}).  \label{est}
\end{equation}%
Here we used the fact that $\widetilde{G}_{n}$ solves the elliptic type problem (\ref{GS_NS}%
) for $f=G(t,Y_{n})$ and assumption (\ref{LG})$_{2}.$

Let us take $t\in \lbrack 0,T]$, the integration over the time interval $%
(0,s)$, $0\leq s\leq \tau _{N}^{n}\wedge t$ of equality (\ref{sh1}) and
estimate (\ref{est}) yield
\begin{align}
\left\Vert Y_{n}(s)\right\Vert _{V}^{2}+\int_{0}^{s}(4\nu \left\Vert
DY_{n}\right\Vert _{2}^{2}\,& +2\nu \gamma \left\Vert Y_{n}\right\Vert
_{L^{2}(\Gamma )}^{2})dr\leq \left\Vert Y_{n}(0)\right\Vert
_{V}^{2}+C(1+\int_{0}^{s}\left\Vert U\right\Vert _{2}^{2}\,dr)  \notag \\
& +\int_{0}^{s}\left\Vert Y_{n}\right\Vert _{V}^{2}\,dr+2\int_{0}^{s}\left(
G\left( r,Y_{n}\right) ,Y_{n}\right) \,dW_{r}.  \label{y2}
\end{align}%
The Burkholder-Davis-Gundy inequality gives
\begin{align*}
\mathbb{E}\sup_{s\in \lbrack 0,\tau _{N}^{n}\wedge t]}\left\vert
\int_{0}^{s}\left( G\left( r,Y_{n}\right) ,Y_{n}\right) \,dW_{r}\right\vert
& \leq \mathbb{E}\left( \int_{0}^{\tau _{N}^{n}\wedge t}\left\vert \left(
G\left( s,Y_{n}\right) ,Y_{n}\right) \right\vert ^{2}\,ds\right) ^{\frac{1}{2%
}} \\
& \leq \mathbb{E}\sup_{s\in \lbrack 0,\tau _{N}^{n}\wedge t]}\left\Vert
Y_{n}(s)\right\Vert _{2}\left( \int_{0}^{\tau _{N}^{n}\wedge t}\left\Vert
G\left( s,Y_{n}\right) \right\Vert _{2}^{2}\,ds\right) ^{\frac{1}{2}} \\
& \leq \varepsilon \,\mathbb{E}\sup_{s\in \lbrack 0,\tau _{N}^{n}\wedge
t]}\Vert Y_{n}(s)\Vert _{2}^{2}+C_{\varepsilon }\mathbb{E}\int_{0}^{\tau
_{N}^{n}\wedge t}(1+\left\Vert Y_{n}\right\Vert _{V}^{2})\,ds.
\end{align*}%
Substituting the last inequality with the chosen $\varepsilon =\frac{1}{2}$
in (\ref{y2}) and considering (\ref{CI}), we derive
\begin{eqnarray*}
\frac{1}{2}\mathbb{E}\sup_{s\in \lbrack 0,\tau _{N}^{n}\wedge t]}\Vert
Y_{n}(s)\Vert _{V}^{2} &+&\mathbb{E}\int_{0}^{\tau _{N}^{n}\wedge t}\left(
4\nu \left\Vert DY_{n}\right\Vert _{2}^{2}\,+2\nu \gamma \left\Vert
Y_{n}\right\Vert _{L^{2}(\Gamma )}^{2}\right) \,ds\leq \mathbb{E}\left\Vert
Y_{0}\right\Vert _{V}^{2} \\
&+&C\mathbb{E}\int_{0}^{t}(1+\left\Vert U\right\Vert _{2}^{2})\,ds+C\mathbb{E%
}\int_{0}^{\tau _{N}^{n}\wedge t}\left\Vert Y_{n}\right\Vert _{V}^{2}\,ds.
\end{eqnarray*}%
Hence, if we denote by $1_{[0,\tau _{N}^{n}]}$ the characteristic function
of the interval $[0,\tau _{N}^{n}],$ the function
\begin{equation*}
f(t)=\mathbb{E}\sup_{s\in \lbrack 0,t]}1_{[0,\tau _{N}^{n}]}\Vert
Y_{n}(s)\Vert _{V}^{2}
\end{equation*}%
fulfills Gronwall%
\'{}%
s type inequality%
\begin{equation*}
\frac{1}{2}f(t)\leq C\int_{0}^{t}f(s)ds+\mathbb{E}\left\Vert
Y_{n}(0)\right\Vert _{V}^{2}+C\mathbb{E}\int_{0}^{t}(1+\left\Vert
U\right\Vert _{2}^{2})\,ds,
\end{equation*}%
which implies
\begin{eqnarray}
\frac{1}{2}\mathbb{E}\sup_{s\in \lbrack 0,\tau _{N}^{n}\wedge t]}\left\Vert
Y_{n}(s)\right\Vert _{V}^{2} &+&\mathbb{E}\int_{0}^{\tau _{N}^{n}\wedge
t}\left( 4\nu \left\Vert DY_{n}\right\Vert _{2}^{2}\,+2\nu \gamma \left\Vert
Y_{n}\right\Vert _{L^{2}(\Gamma )}^{2}\right) \,ds  \notag \\
&\leq &C\left( 1+\mathbb{E}\left\Vert Y_{0}\right\Vert _{V}^{2}+\mathbb{E}%
\Vert U\Vert _{L^{2}(0,t;L^{2})}^{2}\right) .  \label{IMP_1}
\end{eqnarray}


\textit{Step 2. $L^{2}$ \ estimate for $\mathrm{curl}\,$}$\upsilon (Y_{n})$%
\textit{, depending on the stopping times }$\tau _{N}^{n}.$\textit{\ }

The deduction of this estimate is  quite long. Let us first consider the
solutions $\tilde{f}_{n}$ and $\tilde{G}_{n}$ of (\ref{GS_NS}) for $%
f=f(Y_{n})$ and $f=G(t,Y_{n})$, respectively. Then the following relations
hold
\begin{equation}
(\tilde{f}_{n},e_{i})_{V}=(f(Y_{n}),e_{i}),\qquad (\tilde{G}%
_{n},e_{i})_{V}=(G(t,Y_{n}),e_{i}).  \label{gs1}
\end{equation}%
If we use these relations in equality (\ref{y3}), we get
\begin{equation*}
d\left( Y_{n},e_{i}\right) _{V}=\left( \tilde{f}_{n},e_{i}\right)
_{V}\,dt+\left( \tilde{G}_{n},e_{i}\right) _{V}\,dW_{t}.
\end{equation*}%
Multiplying the last identity by $\lambda _{i}$ and using  (\ref{y4}) in the resulting equation yields
\begin{equation*}
d\left( Y_{n},e_{i}\right) _{\widetilde{W}}=(\tilde{f}_{n},e_{i})_{%
\widetilde{W}}\,dt+(\tilde{G}_{n},e_{i})_{\widetilde{W}}\,dW_{t}.
\end{equation*}%
 On the other hand, the It\^{o} formula gives
\begin{equation*}
d\left( Y_{n},e_{i}\right) _{\widetilde{W}}^{2}=2\left( Y_{n},e_{i}\right) _{%
\widetilde{W}}(\tilde{f}_{n},e_{i})_{\widetilde{W}}\,dt+2\left(
Y_{n},e_{i}\right) _{\widetilde{W}}(\tilde{G}_{n},e_{i})_{\widetilde{W}%
}\,dW_{t}+|(\tilde{G}_{n},e_{i})_{\widetilde{W}}|^{2}\,dt.
\end{equation*}%
Multiplying this equality by $\frac{1}{\lambda _{i}}$ and summing over $%
i=1,\dots ,n$, we obtain
\begin{equation*}
d\left\Vert Y_{n}\right\Vert _{\widetilde{W}}^{2}=2(\tilde{f}_{n},Y_{n})_{%
\tilde{W}}\,dt+2(\tilde{G}_{n},Y_{n})_{\widetilde{W}}\,dW_{t}+\sum_{i=1}^{n}%
\frac{1}{\lambda _{i}}|(\tilde{G}_{n},e_{i})_{\widetilde{W}}|^{2}\,dt,
\end{equation*}%
that is
\begin{align*}
d(\Vert \mathrm{curl}\,\upsilon \left( Y_{n}\right) \Vert _{2}^{2}& +\Vert
Y_{n}\Vert _{V}^{2})=2\left( (\mathrm{curl}\,\upsilon (\tilde{f}_{n}),%
\mathrm{curl}\,\upsilon (Y_{n}))+(\tilde{f}_{n},Y_{n})_{V}\right) \,dt \\
& +2\left( (\mathrm{curl}\,\upsilon (\tilde{G}_{n}),\mathrm{curl}\,\upsilon
(Y_{n}))+(\tilde{G}_{n},Y_{n})_{V}\right) \,dW_{t}+\sum_{i=1}^{n}\lambda
_{i}|(G(t,Y_{n}),e_{i})_{V}|^{2}\,dt
\end{align*}%
by the definition of the inner product (\ref{inpr}). The definition of $%
\tilde{f}_{n}$ and $\tilde{G}_{n}$ as solutions of (\ref{GS_NS}) implies
\begin{align*}
d(\Vert \mathrm{curl}\,\upsilon \left( Y_{n}\right) \Vert _{2}^{2}& +\Vert
Y_{n}\Vert _{V}^{2})=2\left( (\mathrm{curl}\,f(Y_{n}),\mathrm{curl}%
\,\upsilon (Y_{n}))+(f(Y_{n}),Y_{n})\right) \,dt \\
& +2\left( (\mathrm{curl}\,G(t,Y_{n}),\mathrm{curl}\,\upsilon
(Y_{n}))+(G(t,Y_{n}),Y_{n})\right) \,dW_{t}+\sum_{i=1}^{n}\lambda
_{i}|(G(t,Y_{n}),e_{i})_{V}|^{2}\,dt,
\end{align*}%
that reduces to
\begin{align}
d\Vert \mathrm{curl}\,\upsilon \left( Y_{n}\right) \Vert _{2}^{2}& =2\left( (%
\mathrm{curl}\,f(Y_{n}),\mathrm{curl}\,\upsilon (Y_{n}))\right) \,dt  \notag
\\
& +2\left( (\mathrm{curl}\,G(t,Y_{n}),\mathrm{curl}\,\upsilon
(Y_{n}))\right) \,dW_{t}+\sum_{i=1}^{n}(\lambda
_{i}-1)|(G(t,Y_{n}),e_{i})|^{2}\,dt,  \label{red}
\end{align}%
taking into account equality (\ref{sh1}).

Since
\begin{equation*}
\mathrm{curl}\left( \mathrm{curl}\left( \upsilon (Y_{n})\right) \times
Y_{n}\right) =\left( Y_{n}\cdot \nabla \right) \mathrm{curl}\,\upsilon
(Y_{n})\qquad \text{and\quad }\left( \left( Y_{n}\cdot \nabla \right)
\mathrm{curl}\,\upsilon (Y_{n}),\mathrm{curl}\,\upsilon (Y_{n})\right) =0,
\end{equation*}%
we have%
\begin{align}
\left( \mathrm{curl}\,f(Y_{n}),\mathrm{curl}\,\upsilon (Y_{n})\right) &
=\left( \nu \,\mathrm{curl}\,\Delta Y_{n}+\mathrm{curl}\,U,\mathrm{curl}%
\,\upsilon (Y_{n})\right)  \notag \\
& =\left( -\frac{\nu }{\alpha }\,\mathrm{curl}\,\upsilon (Y_{n})+\frac{\nu }{%
\alpha }\mathrm{curl}\,Y_{n}+\mathrm{curl}\,U,\mathrm{curl}\,\upsilon
(Y_{n})\right) .  \notag
\end{align}%
Substituting this last relation in (\ref{red}), we derive%
\begin{align}
d\left\Vert \mathrm{curl}\,\upsilon \left( Y_{n}\right) \right\Vert
_{2}^{2}& +\frac{2\nu }{\alpha }\left\Vert \mathrm{curl}\,\upsilon \left(
Y_{n}\right) \right\Vert _{2}^{2}\,dt=2\left( \frac{\nu }{\alpha }\mathrm{%
curl}\,Y_{n}+\mathrm{curl}\,U,\mathrm{curl}\,\upsilon (Y_{n})\right) \,dt
\notag \\
& +2\left( (\mathrm{curl}\,G(t,Y_{n}),\mathrm{curl}\,\upsilon
(Y_{n}))\right) \,dW_{t}+\sum_{i=1}^{n}(\lambda
_{i}-1)|(G(t,Y_{n}),e_{i})|^{2}\,dt.  \label{temps_n}
\end{align}%
Let us take $t\in[0,T]$. By integrating over the time interval $(0,s)$, $%
0\leq s\leq \tau^n _{N}\wedge t$, taking the supremum and the expectation,
we get%
\begin{align}
\mathbb{E}\sup_{s\in \lbrack 0,\tau^n _{N}\wedge t]}\left\Vert \mathrm{curl}%
\,\upsilon \left( Y_{n}(s)\right) \right\Vert _{2}^{2}& +\frac{2\nu }{\alpha
}\mathbb{E}\int_{0}^{\tau^n _{N}\wedge t}\left\Vert \mathrm{curl}\,\upsilon
\left( Y_{n}\right) \right\Vert _{2}^{2}\,ds\leq \mathbb{E}\left\Vert
\mathrm{curl}\,\upsilon \left( Y_{n}(0)\right) \right\Vert _{2}^{2}  \notag
\\
& +2\mathbb{E}\,\int_{0}^{\tau^n _{N}\wedge t}\left\vert \left( \frac{\nu }{%
\alpha }\mathrm{curl}\,Y_{n}+\mathrm{curl}\,U,\mathrm{curl}\,\upsilon
(Y_{n})\right) \right\vert \,ds  \notag \\
& +2\mathbb{E}\,\sup_{s\in \lbrack 0,\tau^n _{N}\wedge t]}\left\vert
\int_{0}^{s}\left( \mathrm{curl}\,G(r,Y_{n}),\mathrm{curl}\,\upsilon
(Y_{n})\right) \,dW_{r}\right\vert  \notag \\
& +\mathbb{E}\,\int_{0}^{\tau^n _{N}\wedge t}\sum_{i=1}^{n}|\lambda
_{i}-1|\, |(G(s,Y_{n}),e_{i})|^{2}\,ds.  \label{4.12}
\end{align}

\bigskip

Moreover
\begin{align*}
& 2\mathbb{E}\,\int_{0}^{\tau^n _{N}\wedge t}\left\vert \left( \frac{\nu }{%
\alpha }\mathrm{curl}\,Y_{n}+\mathrm{curl}\,U,\mathrm{curl}\,\upsilon
(Y_{n})\right) \right\vert \,ds \\
& \leq \mathbb{E}\sup_{s\in \lbrack 0,\tau^n _{N}\wedge t]}\left\Vert
\mathrm{curl}\,\upsilon \left( Y_{n}(s)\right) \right\Vert
_{2}\int_{0}^{\tau _{N}\wedge t}\left( \frac{\nu }{\alpha }\left\Vert
\mathrm{curl}\,Y_{n}\right\Vert _{2}+\left\Vert \mathrm{curl}\,U\right\Vert
_{2}\right) \,ds \\
& \leq \varepsilon \mathbb{E}\sup_{s\in \lbrack 0,\tau^n _{N}\wedge
t]}\left\Vert \mathrm{curl}\,\upsilon \left( Y_{n}(s)\right) \right\Vert
_{2}^{2}+C_{\varepsilon }\mathbb{E}\int_{0}^{\tau^n _{N}\wedge t}(\left\Vert
\mathrm{curl}\,Y_{n}\right\Vert _{2}^{2}+\left\Vert \mathrm{curl}%
\,U\right\Vert _{2}^{2})\,ds,
\end{align*}%
that is
\begin{align}
& 2\mathbb{E}\,\int_{0}^{\tau^n _{N}\wedge t}\left\vert \left( \frac{\nu }{%
\alpha }\mathrm{curl}\,Y_{n}+\mathrm{curl}\,U,\mathrm{curl}\,\upsilon
(Y_{n})\right) \right\vert \,ds  \notag \\
& \leq \varepsilon \mathbb{E}\sup_{s\in \lbrack 0,\tau^n _{N}\wedge
t]}\left\Vert \mathrm{curl}\,\upsilon \left( Y_{n}(s)\right) \right\Vert
_{2}^{2}+C_{\varepsilon }\mathbb{E}\int_{0}^{\tau^n _{N}\wedge t}(\left\Vert
\mathrm{curl}\,Y_{n}\right\Vert _{2}^{2}+\left\Vert \mathrm{curl}%
\,U\right\Vert _{2}^{2})\,ds.  \label{est1}
\end{align}

\bigskip

The Burkholder-Davis-Gundy inequality and estimate (\ref{LG})$_{2}$ imply
\begin{align*}
2\mathbb{E}\sup_{s\in \lbrack 0,\tau^n _{N}\wedge t]}& \left\vert
\int_{0}^{s}\left( \mathrm{curl}\,G\left( r,Y_{n}\right) ,\mathrm{curl}%
\,\upsilon (Y_{n})\right) \,dW_{r}\right\vert \leq 2\mathbb{E}\left(
\int_{0}^{\tau^n _{N}\wedge t}\left\vert \left( \mathrm{curl}\,G\left(
s,Y_{n}\right) ,\mathrm{curl}\,\upsilon (Y_{n})\right) \right\vert
^{2}\,ds\right) ^{\frac{1}{2}} \\
& \leq 2\mathbb{E}\sup_{s\in \lbrack 0,\tau^n _{N}\wedge t]}\left\Vert
\mathrm{curl}\,\upsilon (Y_{n}(s))\right\Vert _{2}\left( \int_{0}^{\tau^n
_{N}\wedge t}\left\Vert \mathrm{curl}\,G\left( s,Y_{n}\right) \right\Vert
_{2}^{2}\,ds\right) ^{\frac{1}{2}} \\
& \leq \varepsilon \,\mathbb{E}\sup_{s\in \lbrack 0,\tau^n _{N}\wedge
t]}\Vert \mathrm{curl}\,\upsilon (Y_{n})\Vert _{2}^{2}+C_{\varepsilon }%
\mathbb{E}\int_{0}^{\tau^n _{N}\wedge t}\left\Vert G\left( s,Y_{n}\right)
\right\Vert _{V}^{2}\,ds \\
& \leq \varepsilon \,\mathbb{E}\sup_{s\in \lbrack 0,\tau^n _{N}\wedge
t]}\Vert \mathrm{curl}\,\upsilon (Y_{n})\Vert _{2}^{2}+C_{\varepsilon }%
\mathbb{E}\int_{0}^{\tau^n _{N}\wedge t}\left( 1+\left\Vert
\,Y_{n}(s)\right\Vert _{V}^{2}\right) \,ds,
\end{align*}%
that is
\begin{align}
2\mathbb{E}\sup_{s\in \lbrack 0,\tau^n _{N}\wedge t]}& \left\vert
\int_{0}^{s}\left( \mathrm{curl}\,G\left( r,Y_{n}\right) ,\mathrm{curl}%
\,\upsilon (Y_{n})\right) \,dW_{r}\right\vert  \notag \\
& \leq \varepsilon \,\mathbb{E}\sup_{s\in \lbrack 0,\tau^n _{N}\wedge
t]}\Vert \mathrm{curl}\,\upsilon (Y_{n})\Vert _{2}^{2}+C_{\varepsilon }%
\mathbb{E}\int_{0}^{\tau^n _{N}\wedge t}\left( 1+\left\Vert
\,Y_{n}\right\Vert _{V}^{2}\right) \,ds.  \label{est2}
\end{align}

\bigskip

Substituting (\ref{est1})-(\ref{est2}) in (\ref{4.12}) and choosing $%
\varepsilon =\frac{1}{4}$, we obtain
\begin{align}
\frac{1}{2}\mathbb{E}\sup_{s\in \lbrack 0,\tau _{N}^{n}\wedge t]}\left\Vert
\mathrm{curl}\,\upsilon \left( Y_{n}(s)\right) \right\Vert _{2}^{2}& +\frac{%
2\nu }{\alpha }\mathbb{E}\int_{0}^{\tau _{N}^{n}\wedge t}\left\Vert \mathrm{%
curl}\,\upsilon \left( Y_{n}\right) \right\Vert _{2}^{2}\,ds\leq \mathbb{E}%
\left\Vert \mathrm{curl}\,\upsilon \left( Y_{n}(0)\right) \right\Vert
_{2}^{2}  \notag \\
& +C\mathbb{E}\int_{0}^{t}\left\Vert \mathrm{curl}\,U\right\Vert
_{2}^{2}\,ds+C\mathbb{E}\int_{0}^{\tau _{N}^{n}\wedge t}\left( 1+\left\Vert
Y_{n}(s)\right\Vert _{V}^{2}\right) \,ds.  \label{IMP_4}
\end{align}

\bigskip

\textit{Step 3. The limit transition, as }$N\rightarrow \infty ,$\textit{\ \
in estimates (\ref{IMP_1}), (\ref{IMP_4}).}\bigskip

Since%
\begin{equation*}
\mathbb{E}\left\Vert \mathrm{curl}\,\upsilon \left( Y_{n}(0)\right)
\right\Vert _{2}^{2}\leq C\mathbb{E}\left\Vert Y_{n}(0)\right\Vert
_{H^{3}}^{2}\leq C\mathbb{E}\left\Vert Y_{0}\right\Vert _{H^{3}}^{2}\leq C
\end{equation*}%
and
\begin{equation*}
\mathbb{E}\int_{0}^{\tau _{N}^{n}\wedge t}\left( 1+\left\Vert
Y_{n}(s)\right\Vert _{V}^{2}\right) \,ds\leq C
\end{equation*}%
by (\ref{CI}) and (\ref{IMP_1}), we obtain
\begin{equation*}
\mathbb{E}\sup_{s\in \lbrack 0,\tau _{N}^{n}\wedge t]}\left\Vert \mathrm{curl%
}\,\upsilon \left( Y_{n}(s)\right) \right\Vert _{2}^{2}\leq C.
\end{equation*}%
Therefore estimates (\ref{y2_sigma}), (\ref{IMP_1}) imply
\begin{equation*}
\mathbb{E}\sup_{s\in \lbrack 0,\tau _{N}^{n}\wedge T]}\left\Vert
Y_{n}(s)\right\Vert _{H^{3}}^{2}\leq C,
\end{equation*}%
where $C$ is a constant independent of $N$ and $n$. Let us fix $n\in \mathbb{%
N}$, writing
\begin{align}
\mathbb{E}\sup_{s\in \lbrack 0,\tau _{N}^{n}\wedge T]}\left\Vert
Y_{n}(s)\right\Vert _{H^{3}}^{2}& =\mathbb{E}\left( \sup_{s\in \lbrack
0,\tau _{N}^{n}\wedge T]}1_{\{\tau _{N}^{n}<T\}}\ \left\Vert
Y_{n}(s)\right\Vert _{H^{3}}^{2}\right)  \notag \\
& +\mathbb{E}\left( \sup_{s\in \lbrack 0,\tau _{N}^{n}\wedge T]}1_{\{\tau
_{N}^{n}\geq T\}}\ \left\Vert Y_{n}(s)\right\Vert _{H^{3}}^{2}\right)  \notag
\\
& \geq \mathbb{E}\left( \max_{s\in \lbrack 0,\tau _{N}^{n}]}1_{\{\tau
_{N}^{n}<T\}}\ \left\Vert Y_{n}(s)\right\Vert _{H^{3}}^{2}\right)
\geq N^{2}P\left( \tau _{N}^{n}<T\right) ,  \label{fc}
\end{align}%
we deduce that $P\left( \tau _{N}^{n}<T\right) \rightarrow 0$, as $%
N\rightarrow \infty .$ This means that $\tau _{N}^{n}\rightarrow T$ in
probability, as $N\rightarrow \infty $. Then there exists a subsequence $%
\{\tau _{N_{k}}^{n}\}$ of $\{\tau _{N}^{n}\}$ (that may depend on $n$) such
that
\begin{equation*}
\tau _{N_{k}}^{n}(\omega )\rightarrow T\text{\qquad for a. e. \ }\omega \in
\Omega \quad \text{as }k\rightarrow \infty .
\end{equation*}%
Since $\tau _{N_{k}}^{n}\leq T_{n}\leq T$, we deduce that $T_{n}=T$, so $%
Y_{n}$ is a global-in-time solution of the stochastic differential equation (%
\ref{y1}). On the other hand, the sequence $\left\{ \tau _{N}^{n}\right\} $
of the stopping times is monotone on $N$ for each fixed $n$, then we can
apply the monotone convergence theorem in order to pass to the limit in
inequalities (\ref{IMP_1}) and (\ref{IMP_4}) as $N\rightarrow \infty $, \
deducing estimates (\ref{ineq1}) and (\ref{ineq222}).

\bigskip

\textit{Step 4. Estimate in the space }$\widetilde{W}$\ \textit{\ for }$Y_n$%
\textit{. } Substituting estimate (\ref{ineq1}) in (\ref{ineq222}) and using
Lemma \ref{non_lin}, we immediately derive the main estimate \ (\ref{ineq2})
of this lemma.\bigskip\ $\hfill \hfill \blacksquare $

\bigskip

In the next lemma, assuming a better integrability for the data $U$, $Y_{0},$
we improve the integrability properties for the solution $Y_{n}$ of problem (%
\ref{y1}).

\begin{lemma}
\label{existence_state1} Assume that%
\begin{equation*}
U\in L^{p}(\Omega \times (0,T);H),\qquad Y_{0}\in L^{p}(\Omega ,V)\qquad
\text{for some\ }\ 4\leq p<\infty .
\end{equation*}%
Then the solution $Y_{n}$ of problem $(\ref{y1})$ belongs to $L^{p}(\Omega
,L^{\infty }(0,T;V))$ and verifies the estimate%
\begin{equation}
\mathbb{E}\sup_{s\in \lbrack 0,t]}\left\Vert Y_{n}(t)\right\Vert
_{V}^{p}\leq C\mathbb{E}\left\Vert Y_{0}\right\Vert _{V}^{p}+C\,(1+\mathbb{E}%
\int_{0}^{t}\left\Vert U\right\Vert _{2}^{p}\,ds),  \label{lp1}
\end{equation}%
where $C$ is a positive constant independent of $n.$
\end{lemma}

\textbf{Proof.} For each $n\in \mathbb{N},$ let us define the suitable
sequence $\{\tau _{N}^{n}\}_{N\in \mathbb{N}}$ of the stopping times
\begin{equation*}
\tau _{N}^{n}=\inf \{t\geq 0:\Vert Y_{n}(t)\Vert _{V}\geq N\}\wedge T.
\end{equation*}%
Applying the It\^{o} formula for the function $\theta (x)=x^{q},\,q\geq 1,$
to process (\ref{sh1}), we have%
\begin{align*}
d\left\Vert Y_{n}\right\Vert _{V}^{2q}& =q\left\Vert Y_{n}\right\Vert
_{V}^{2q-2}\left[ -\left( 4\nu \left\Vert DY_{n}\right\Vert _{2}^{2}\,+2\nu
\gamma \left\Vert Y_{n}\right\Vert _{L^{2}(\Gamma )}^{2}\right) +2\left(
U,Y_{n}\right) +\sum_{i=1}^{n}|\left( G\left( t,Y_{n}\right) ,e_{i}\right)|
^{2}\right] \,dt \\
& +2q\left\Vert Y_{n}\right\Vert _{V}^{2q-2}\left( G(t,Y_{n}),Y_{n}\right)
\,dW_{t}+2q(q-1)\left\Vert Y_{n}\right\Vert _{V}^{2q-4} |\left(
G(t,Y_{n}),Y_{n}\right)| ^{2}\,dt.
\end{align*}%
Let us take $t\in \lbrack 0,T]$. Integrating over the time interval $[0,s]$,
$0\leq s\leq \tau _{N}^{n}\wedge t$, we obtain
\begin{align}
\left\Vert Y_{n}(s)\right\Vert _{V}^{2q}& \leq \left\Vert
Y_{n}(0)\right\Vert _{V}^{2q}+q\int_{0}^{s}\left\Vert Y_{n}\right\Vert
_{V}^{2q-2}\left\vert 2\left( U,Y_{n}\right) +\sum_{i=1}^{n}|\left( G\left(
r,Y_{n}\right) ,e_{i}\right)| ^{2}\right\vert \,dr  \notag \\
& +2q\left\vert \int_{0}^{s}\left\Vert Y_{n}\right\Vert _{V}^{2q-2}\left(
G(r,Y_{n}),Y_{n}\right) \,dW_{r}\right\vert  \notag \\
& +2q(q-1)\int_{0}^{s}\left\Vert Y_{n}\right\Vert _{V}^{2q-4}|\left(
G(r,Y_{n}),Y_{n}\right)| ^{2}\,dr.  \label{r1}
\end{align}

From estimate (\ref{est}) we have
\begin{equation*}
\sum_{i=1}^{n}|\left( G\left( t,Y_{n}\right) ,e_{i}\right)| ^{2}\leq
C(1+||Y_{n}||_{V}^{2}).
\end{equation*}%
Taking the supremum on $s\in \lbrack 0,\tau _{N}^{n}\wedge t]$, \ the
expectation in (\ref{r1}), applying Burkholder-Davis-Gundy%
\'{}%
s and Young%
\'{}%
s inequalities, and proceeding analogously to (\ref{y2}), we obtain
\begin{align*}
\mathbb{E}\sup_{s\in \lbrack 0,\tau _{N}^{n}\wedge t]}\left\Vert
Y_{n}(s)\right\Vert _{V}^{2q}& \leq \mathbb{E}\left\Vert Y_{n}(0)\right\Vert
_{V}^{2q}+C_{q}\,\mathbb{E}\int_{0}^{t}\left\Vert U\right\Vert _{2}^{2q}\,ds
\\
& +\frac{1}{2}\mathbb{E}\sup_{s\in \lbrack 0,\tau _{N}^{n}\wedge
t]}\left\Vert Y_{n}(s)\right\Vert _{V}^{2q}+C\,_{q}(1+\,\mathbb{E}%
\int_{0}^{\tau _{N}\wedge t}\left\Vert Y_{n}\right\Vert _{V}^{2q}\,ds).
\end{align*}%
Using Gronwall%
\'{}%
s inequality, we deduce
\begin{equation}
\mathbb{E}\sup_{s\in \lbrack 0,\tau _{N}^{n}\wedge t]}\left\Vert
Y_{n}(s)\right\Vert _{V}^{2q}\leq C\mathbb{E}\left\Vert Y_{0}\right\Vert
_{V}^{2q}+C\,(1+\mathbb{E}\int_{0}^{t}\left\Vert U\right\Vert _{2}^{2q}\,ds)
\label{IMP_7}
\end{equation}%
for any $q\geq 1$ and $t\in \lbrack 0,T]$. Using the fact that
\begin{equation*}
\mathbb{E}\sup_{s\in \lbrack 0,\tau _{N}^{n}\wedge T]}\left\Vert
Y_{n}(s)\right\Vert _{V}^{2q}\leq C
\end{equation*}%
with $C$ independent of $n$ and $N$, we may reasoning as in the proof of
Lemma \ref{existence_state}, in order to verify that for each $n$, $\tau
_{N}^{n}\rightarrow T$ in probability, as $N\rightarrow \infty $. Then,
there exists a subsequence $\{\tau _{N_{k}}^{n}\}$ of $\{\tau _{N}^{n}\}$
(that may depend on $n$) such that $\tau _{N_{k}}^{n}\rightarrow T$ \ for a.
e. $\omega \in \Omega $, as $k\rightarrow \infty $. Now,  let us consider  $q=\frac{p}{2}$. Using the monotone
convergence theorem, we pass to the limit in (\ref{IMP_7}) as $k\rightarrow
\infty $, deriving estimate (\ref{lp1}).$\hfill \hfill \blacksquare $

\bigskip

\textbf{Proof of Theorem \ref{the_1}.} \vspace{2mm} \textbf{Existence.} The
proof is split into four steps. \vspace{2mm}\newline
\textit{Step 1. \ Estimates and convergences, related with the projection
operator.}

Let $P_{n}:\widetilde{W}\rightarrow W_{n}$ be the orthogonal projection
defined by%
\begin{equation*}
P_{n}y=\sum_{j=1}^{n}\widetilde{c}_{j}\widetilde{e}_{j}\quad \text{ with }%
\quad \text{ }\widetilde{c}_{j}=\left( y,\widetilde{e}_{j}\right) _{%
\widetilde{W}},\qquad \forall y\in \widetilde{W},
\end{equation*}%
where $\{\widetilde{e}_{j}=\frac{1}{\sqrt{\lambda _{j}}}e_{j}\}_{j=1}^{%
\infty }$ is the orthonormal basis of $\ \widetilde{W}.$ \ \bigskip It is
easy to check that
\begin{equation*}
P_{n}y=\sum_{j=1}^{n}c_{j}e_{j}\quad \text{ with }\quad \text{ }c_{j}=\left(
y,e_{j}\right) _{V},\qquad \forall y\in \widetilde{W}.
\end{equation*}

By Parseval's identity we have that
\begin{equation*}
||P_{n}y||_{V}\leq ||y||_{V},\text{ }\quad \forall y\in V,
\end{equation*}
\begin{equation*}
||P_{n}y||_{\widetilde{W}}\leq ||y||_{\widetilde{W}}\qquad \mbox{ and
}\qquad P_{n}y\longrightarrow y\qquad \mbox{ strongly in
}\ \widetilde{W},\text{ }\quad \forall y\in \widetilde{W}.
\end{equation*}

Considering an arbitrary $Z\in L^{2}(\Omega \times (0,T);\widetilde{W})),$
we have%
\begin{equation*}
||P_{n}Z||_{\widetilde{W}}\leq ||Z||_{\widetilde{W}}\qquad \mbox{ and
}\qquad P_{n}Z(\omega ,t)\rightarrow Z(\omega ,t)\qquad \mbox{ strongly in
}\ \widetilde{W},
\end{equation*}%
which are valid for $P$-a. e. $\omega \in \Omega $ and a.e. $t\in (0,T).$
Hence Lebesgue's dominated convergence theorem and the inequality
\begin{equation*}
||Z||_{V}\leq C||Z||_{\widetilde{W}}\qquad \text{for any }Z\in \widetilde{W}
\end{equation*}%
imply%
\begin{eqnarray}
P_{n}Z &\longrightarrow &Z\qquad \mbox{ strongly in
}\ L^{2}(\Omega \times (0,T);\widetilde{W})),  \notag \\
P_{n}Z &\longrightarrow &Z\qquad \mbox{ strongly in
}\ L^{2}(\Omega \times (0,T);V)).  \label{c02}
\end{eqnarray}%
\bigskip \textit{Step 2.} \textit{Passing to the limit in the weak sense.}

We have
\begin{equation}
\mathbb{E}\sup_{t\in \lbrack 0,T]}\left\Vert Y_{n}(t)\right\Vert _{%
\widetilde{W}}^{2}\leq C,\qquad \mathbb{E}\sup_{t\in \lbrack 0,T]}\left\Vert
Y_{n}(t)\right\Vert _{V}^{p}\leq C  \label{as}
\end{equation}%
for some constants $C$ independent of the index $n$, by estimates (\ref%
{ineq2}) and (\ref{lp1}). Therefore there exists a suitable subsequence $%
Y_{n}$, which is indexed by the same index $n$, for simplicity of notations,
such that
\begin{eqnarray}
Y_{n} &\rightharpoonup &Y\qquad \mbox{ *-weakly in
}\ L^{2}(\Omega ,L^{\infty }(0,T;\widetilde{W})),  \notag \\
Y_{n} &\rightharpoonup &Y\qquad \mbox{ *-weakly in
}\ L^{p}(\Omega ,L^{\infty }(0,T;V)).  \label{c1}
\end{eqnarray}%
Moreover, we have
\begin{eqnarray}
P_{n}Y &\longrightarrow &Y\qquad \mbox{ strongly in
}\ L^{2}(\Omega \times (0,T);\widetilde{W})),  \notag \\
P_{n}Y &\longrightarrow &Y\qquad \mbox{ strongly in
}\ L^{2}(\Omega \times (0,T);V)).  \label{c02Y}
\end{eqnarray}%
Let us introduce the operator $B:\widetilde{W}\times V\rightarrow \widetilde{%
W}^{\ast },$ defined as
\begin{equation*}
\left( B\left( y,z\right) ,\phi \right) =\left( \mathrm{curl}\,\upsilon
(y)\times z,\phi \right) \qquad \mbox{ for any
}\ y,\phi \in \widetilde{W}\qquad \mbox{ and
}\ z\in V,
\end{equation*}%
and state some useful properties. Relation (\ref{bb1}) gives
\begin{equation}
\left( B\left( y,z\right) ,\phi \right) =-\left( B\left( y,\phi \right)
,z\right) ,\quad \quad \left( B\left( y,z\right) ,z\right) =0,  \label{antis}
\end{equation}%
and (\ref{rm2_lin1}), (\ref{rm2_lin3}) yield
\begin{equation}
\Vert B\left( y,z\right) \Vert _{\widetilde{W}^{\ast }}\leq C\left\Vert
z\right\Vert _{V}\left\Vert y\right\Vert _{\widetilde{W}},  \label{proB_2}
\end{equation}%
\begin{equation}
\Vert B\left( y,z\right) \Vert _{\widetilde{W}^{\ast }}\leq C\left\Vert
y\right\Vert _{V}\left\Vert z\right\Vert _{\widetilde{W}}.  \label{proB_3}
\end{equation}%
From (\ref{rm2_lin}) there exists a fixed constant $C_{1}$ such that
\begin{equation}
\Vert B\left( y,y\right) \Vert _{\widetilde{W}^{\ast }}\leq C_{1}\left\Vert
y\right\Vert _{V}^{2},  \label{proB}
\end{equation}%
then
\begin{equation}
\left\Vert B(y,y)\right\Vert _{L^{2}(\Omega \times (0,T);\widetilde{W}^{\ast
}))}\leq C_{1}\left\Vert y\right\Vert _{L^{4}(\Omega ,L^{4}(0,T;V))}^{2}.
\label{l}
\end{equation}%
On the other hand, taking into account (\ref{LG}), (\ref{c1}), there exist operators $%
B^{\ast }(t)$ and $G^{\ast }(t)$, such that%
\begin{eqnarray}
B(Y_{n},Y_{n}) &\rightharpoonup &B^{\ast }(t)\qquad \mbox{ weakly in
}\ L^{2}(\Omega \times (0,T);\widetilde{W}^{\ast }),  \notag \\
G(t,Y_{n}) &\rightharpoonup &G^{\ast }(t)\qquad \mbox{ weakly in
}\ L^{2}(\Omega \times (0,T);V^m).  \label{c01}
\end{eqnarray}%
Passing on the limit $n\rightarrow \infty $ in equation (\ref{y1}), we
derive that the limit function $Y$ satisfies \ the stochastic differential
equation
\begin{equation}
d\left( \upsilon \left( Y\right) ,\phi \right) =\left[ \left( \nu \Delta
Y+U,\phi \right) -\langle B^{\ast }(t),\phi \rangle \right] \,dt\vspace{2mm}%
+\left( G^{\ast }(t),\phi \right) \,dW_{t},\qquad \forall \phi \in
\widetilde{W}.  \label{y12}
\end{equation}

\bigskip

\textit{Step 3.} \textit{Deduction of strong convergences, as }$n\rightarrow
\infty $\textit{, depending on the stopping times }$\tau _{M}.$\textit{\ \ }

In order to prove that the limit process $Y$ satisfy equation (\ref%
{var_form_state}), we will adapt the methods in \cite{B99} (see also \cite%
{RS12}). Let us introduce a sequence $(\tau _{M})$, $M\in \mathbb{N}$, of
stopping times defined by
\begin{equation*}
\tau _{M}(\omega )=\inf \{t\geq 0:\Vert Y(t)\Vert _{\widetilde{W}}(\omega
)\geq M\}\wedge T,\qquad \omega \in \Omega .
\end{equation*}%
Taking the difference of (\ref{y1}) and (\ref{y12}), we deduce
\begin{eqnarray}
d\left( P_{n}Y-Y_{n},e_{i}\right) _{V} &=&\left[ \left( \nu \Delta
(Y-Y_{n}),e_{i}\right) +\langle B(Y_{n},Y_{n})-B^{\ast }(t),e_{i}\rangle %
\right] \,dt  \notag \\
&&\vspace{2mm}-\left( G(t,Y_{n})-G^{\ast }(t),e_{i}\right) \,dW_{t},
\label{diff}
\end{eqnarray}%
which is valid for any $e_{i}\in W_{n},$ $i=1,...,n.$

By applying It\^{o}'s formula, equation (\ref{diff}) gives
\begin{eqnarray*}
d(P_{n}Y-Y_{n},e_{i})_{V}^{2} &=&2\left( P_{n}Y-Y_{n},e_{i}\right) _{V}\left[
\left( \nu \Delta (Y-Y_{n}),e_{i}\right) +\langle B(Y_{n},Y_{n})-B^{\ast
}(t),e_{i}\rangle \right] \,dt\vspace{2mm} \\
&-&2\left( P_{n}Y-Y_{n},e_{i}\right) _{V}\left( G(t,Y_{n})-G^{\ast
}(t),e_{i}\right) \,dW_{t}+|\left( G(t,Y_{n})-G^{\ast }(t),e_{i}\right)|
^{2}\,dt,
\end{eqnarray*}%
and summing over the index $i$ from $1$ to $n,$ \ we derive
\begin{align}
&d\left( ||P_{n}Y-Y_{n}||_{V}^{2}\right) +\left( 4\nu
||D(P_{n}Y-Y_{n})||_{2}^{2}+2\nu \gamma \left\Vert P_{n}Y-Y_{n}\right\Vert
_{L^{2}(\Gamma )}^{2}\right) dt  \notag \\
& =2\nu (\Delta (P_{n}Y-Y),P_{n}Y-Y_{n})\,dt  \notag \\
& +2\langle B(Y_{n},Y_{n})-B^{\ast }(t),P_{n}Y-Y_{n}\rangle \,dt  \notag \\
& +\sum_{i=1}^{n}|\left( G(t,Y_{n})-G^{\ast }(t),e_{i}\right)|
^{2}\,dt-2\left( G(t,Y_{n})-G^{\ast }(t),P_{n}Y-Y_{n}\right) \,dW_{t}.
\label{y13}
\end{align}%
Let us notice that
\begin{align}
\langle B(Y_{n},Y_{n})& -B^{\ast }(t),P_{n}Y-Y_{n}\rangle =\langle
B(Y_{n},Y_{n})-B(P_{n}Y,P_{n}Y),P_{n}Y-Y_{n}\rangle  \notag \\
& +\langle B(P_{n}Y,P_{n}Y)-B(Y,Y),P_{n}Y-Y_{n}\rangle +\langle
B(Y,Y)-B^{\ast }(t),P_{n}Y-Y_{n}\rangle  \notag \\
& =I_{1}+I_{2}+I_{3}.  \label{oi}
\end{align}%
Using (\ref{antis}), we derive
\begin{align*}
I_{1}& =\langle B(Y_{n},Y_{n})-B(P_{n}Y,P_{n}Y),P_{n}Y-Y_{n}\rangle \\
& =\langle
B(Y_{n},Y_{n})-B(Y_{n},P_{n}Y)+B(Y_{n},P_{n}Y)-B(P_{n}Y,P_{n}Y),P_{n}Y-Y_{n}%
\rangle \\
& =-\langle B(Y_{n},P_{n}Y-Y_{n}),P_{n}Y-Y_{n}\rangle -\langle
B(P_{n}Y-Y_{n},P_{n}Y),P_{n}Y-Y_{n}\rangle \\
& =\langle B(P_{n}Y-Y_{n},P_{n}Y-Y_{n}),P_{n}Y\rangle ,
\end{align*}%
which along with (\ref{proB}) implies
\begin{equation}
|I_{1}|\leq C_{1}\left\Vert Y\right\Vert _{\widetilde{W}}\left\Vert
P_{n}Y-Y_{n}\right\Vert _{V}^{2}.  \label{BB1}
\end{equation}%
For the term $I_{2}$, we have
\begin{align*}
|I_{2}|& =\left\vert \langle B(P_{n}Y,P_{n}Y)-B(Y,Y),P_{n}Y-Y_{n}\rangle
\right\vert \\
& \leq \left\Vert B(P_{n}Y,P_{n}Y)-B(Y,Y)\right\Vert _{\widetilde{W}^{\ast
}}\left\Vert P_{n}Y-Y_{n}\right\Vert _{\widetilde{W}},
\end{align*}%
and for every $\phi \in W$, it follows from (\ref{proB_2}) and (\ref{proB_3})
that
\begin{align*}
\Vert B(P_{n}Y,P_{n}Y)-B(Y,Y)\Vert _{\widetilde{W}^{\ast }}& \leq \Vert
B(P_{n}Y-Y,P_{n}Y)\Vert _{\widetilde{W}^{\ast }}+\Vert B(Y,P_{n}Y-Y\Vert _{%
\widetilde{W}^{\ast }} \\
& \leq C\Vert Y\Vert _{\widetilde{W}}\Vert P_{n}Y-Y\Vert _{V}.
\end{align*}%
and consequently, we obtain
\begin{equation}
|I_{2}|\leq C\left\Vert Y\right\Vert _{\widetilde{W}}\left\Vert
P_{n}Y-Y\right\Vert _{V}\left\Vert P_{n}Y-Y_{n}\right\Vert _{\widetilde{W}}.
\label{BB22}
\end{equation}

On the other hand, denoting by $\widetilde{G}_{n}$, $\widetilde{G}$ and $%
\widetilde{G}^{\ast }$ the solutions of the Stokes system (\ref{GS_NS}) for $%
f=G(t,Y_{n})$, $f=G(t,Y)$ and $f=G^{\ast }(t)$, respectively, we have
\begin{equation*}
\left( G(t,Y_{n})-G^{\ast }(t),e_{i}\right) =(\widetilde{G}_{n}-\widetilde{G}%
^{\ast },e_{i})_{V},\qquad i=1,2,\dots ,n.
\end{equation*}%
Then
\begin{equation*}
\sum_{i=1}^{n}|\left( G(t,Y_{n})-G^{\ast }(t),e_{i}\right) |^{2}=\Vert P_{n}%
\widetilde{G}_{n}-P_{n}\widetilde{G}^{\ast }\Vert _{V}^{2}.
\end{equation*}%
The standard relation $x^{2}=(x-y)^{2}-y^{2}+2xy\;$ allows to write
\begin{align*}
\Vert P_{n}\widetilde{G}_{n}-P_{n}\widetilde{G}^{\ast }\Vert _{V}^{2}&
=\Vert P_{n}\widetilde{G}_{n}-P_{n}\widetilde{G}\Vert _{V}^{2}-\Vert P_{n}%
\widetilde{G}-P_{n}\widetilde{G}^{\ast }\Vert _{V}^{2} \\
& +2(P_{n}\widetilde{G}_{n}-P_{n}\widetilde{G}^{\ast },P_{n}\widetilde{G}%
-P_{n}\widetilde{G}^{\ast })_{V}.
\end{align*}%
From the properties of the solutions of the Stokes system (\ref{GS_NS}) and (%
\ref{LG}), we have
\begin{equation*}
\Vert P_{n}\widetilde{G}_{n}-P_{n}\widetilde{G}\Vert _{V}^{2}\leq \Vert
\widetilde{G}_{n}-\widetilde{G}\Vert _{V}^{2}\leq \Vert
G(t,Y_{n})-G(t,Y)\Vert _{L^{2}}^{2}\leq K\left\Vert Y_{n}-Y\right\Vert
_{V}^{2},
\end{equation*}%
then, for the fixed constant $C_{2}=2K$, we have
\begin{align}
\Vert P_{n}\widetilde{G}_{n}-P_{n}\widetilde{G}^{\ast }\Vert _{V}^{2}& \leq
K\left\Vert Y_{n}-Y\right\Vert _{V}^{2}-\Vert P_{n}\widetilde{G}-P_{n}%
\widetilde{G}^{\ast }\Vert _{V}^{2}  \notag \\
& +2(P_{n}\widetilde{G}_{n}-P_{n}\widetilde{G}^{\ast },P_{n}\widetilde{G}%
-P_{n}\widetilde{G}^{\ast })_{V}  \notag \\
& \leq C_{2}\left\Vert Y_{n}-P_{n}Y\right\Vert _{V}^{2}+C\left\Vert
P_{n}Y-Y\right\Vert _{V}^{2}-\Vert P_{n}\widetilde{G}-P_{n}\widetilde{G}%
^{\ast }\Vert _{V}^{2}  \notag \\
& +2(P_{n}\widetilde{G}_{n}-P_{n}\widetilde{G}^{\ast },P_{n}\widetilde{G}%
-P_{n}\widetilde{G}^{\ast })_{V}.  \label{GG1}
\end{align}%
The positive constants $C_{1}$ and $C_{2}$ in (\ref{BB1}) and (\ref{GG1}),
are independent of $n$ and may depend on the data: the domain $\mathcal{O},$
\ the regularity of $\Gamma ,$ the physical constants $\nu $, $\alpha $, $%
\gamma $, $K.$

Let us notice that from the convergence results (\ref{c02})-(\ref{c02Y}), (%
\ref{c01}), we can guess that by passing to the limit in equation (\ref{y13}%
), in a suitable way, as $n\rightarrow \infty ,$ all terms containing $%
P_{n}Y-Y$ will vanish on the right hand side of equality (\ref{y13}),
according to relations (\ref{c02}), (\ref{BB22}) and (\ref{GG1}). But the
terms with $Y_{n}-P_{n}Y$ \ will remain. Fortunately, these terms can be
eliminated by the introduction of the auxiliary function
\begin{equation*}
\xi (t)=e^{-C_{2}t-2C_{1}\int_{0}^{t}\left\Vert Y\right\Vert _{\widetilde{W}%
}ds}.
\end{equation*}%
Now, applying It\^{o}'s formula in equality (\ref{y13}), we get
\begin{align}
& d\left( \xi (t)||P_{n}Y-Y_{n}||_{V}^{2}\right) +\xi (t)\left( 4\nu
||D(P_{n}Y-Y_{n})||_{2}^{2}+2\nu \gamma \left\Vert P_{n}Y-Y_{n}\right\Vert _{L^{2}(\Gamma )}^{2}
\right) dt\notag\\
&=2\nu \xi (t)(\Delta (P_{n}Y-Y),P_{n}Y-Y_{n})\,dt  \notag \\
& +2\xi (t)\langle B(Y_{n},Y_{n})-B^{\ast }(t),P_{n}Y-Y_{n}\rangle \,dt+\xi
(t)\sum_{i=1}^{n}|\left( G(t,Y_{n})-G^{\ast }(t),e_{i}\right)| ^{2}\,dt  \notag
\\
& -2\xi (t)\left( G(t,Y_{n})-G^{\ast }(t),P_{n}Y-Y_{n}\right) \,dW_{t}
\notag \\
& -C_{2}\xi (t)||P_{n}Y-Y_{n}||_{V}^{2}\,dt-2C_{1}\xi (t)\left\Vert
Y\right\Vert _{\widetilde{W}}||P_{n}Y-Y_{n}||_{V}^{2}\,dt.  \label{y14}
\end{align}
Integrating it over the time interval $(0,\tau _{M}(\omega )),$ taking the
expectation and applying estimates (\ref{oi}), (\ref{BB1}), (\ref{GG1}), we
deduce
\begin{align*}
& \mathbb{E}\left( \xi (\tau _{M})||P_{n}Y(\tau _{M})-Y_{n}(\tau
_{M})||_{V}^{2}\right) +\mathbb{E}\int_{0}^{\tau _{M}}\xi (s)\Vert P_{n}%
\widetilde{G}-P_{n}\widetilde{G}^{\ast }\Vert _{V}^{2}ds \\
& +\mathbb{E}\int_{0}^{\tau _{M}}\xi (s)\left( 4\nu
||D(P_{n}Y-Y_{n})||_{2}^{2}+2\nu \gamma \left\Vert P_{n}Y-Y_{n}\right\Vert
_{L^{2}(\Gamma )}^{2}\right) ds \\
& \leq 2\nu \mathbb{E}\int_{0}^{\tau _{M}}\xi (s)(\Delta
(P_{n}Y-Y),P_{n}Y-Y_{n})\,ds \\
& +2\mathbb{E}\int_{0}^{\tau _{M}}\xi (s)I_{2}\,ds+2\mathbb{E}\int_{0}^{\tau
_{M}}~\xi (s)I_{3}\,ds \\
& +\mathbb{E}\int_{0}^{\tau _{M}}~\xi (s)\left[ C\left\Vert
P_{n}Y-Y\right\Vert _{V}^{2}+2(P_{n}\widetilde{G}_{n}-P_{n}\widetilde{G}%
^{\ast },P_{n}\widetilde{G}-P_{n}\widetilde{G}^{\ast })_{V}\right] \,ds \\
& =J_{1}+J_{2}+J_{3}+J_{4}.
\end{align*}%
In what follows we show that, for each $M\in \mathbb{N}$, the right hand
side of this inequality goes to zero, as $n\rightarrow \infty $.

Using (\ref{as})-(\ref{c1}) and the properties of the projection $P_{n},$ we
have
\begin{align*}
|J_{1}|& =\bigl|2\nu \mathbb{E}\int_{0}^{T}\xi (s)(1_{[0,\tau
_{M}]}(s)\Delta (P_{n}Y-Y),P_{n}Y-Y_{n})\,ds\bigr| \\
& \leq C\Vert P_{n}Y-Y\Vert _{L^{2}(\Omega \times (0,T);H^{2})}\Vert
P_{n}Y-Y_{n}\Vert _{L^{2}(\Omega \times (0,T);H^{2})} \\
& \leq C\Vert P_{n}Y-Y\Vert _{L^{2}(\Omega \times (0,T);H^{2})}\left( \Vert
Y\Vert _{L^{2}(\Omega \times (0,T);H^{2})}+\Vert Y_{n}\Vert _{L^{2}(\Omega
\times (0,T);H^{2})}\right) \\
& \leq C\Vert P_{n}Y-Y\Vert _{L^{2}(\Omega \times (0,T);H^{2})},
\end{align*}%
which goes to zero, as $n\rightarrow \infty $, by (\ref{c02Y}). Taking into
account estimates (\ref{ineq2}), (\ref{BB22}) and knowing that $1_{[0,\tau
_{M}]}(s)\left\Vert Y(s)\right\Vert _{\widetilde{W}}\leq M$, P - a. e. in $%
\Omega $, we deduce that
\begin{align*}
|J_{2}|& \leq 2\mathbb{E}\left\vert \int_{0}^{\tau _{M}}\xi
(s)I_{2}\,ds\right\vert \\
& \leq 2\mathbb{E}\int_{0}^{T}\xi (s)1_{[0,\tau _{M}]}(s)\left\Vert
Y\right\Vert _{\widetilde{W}}\left\Vert P_{n}Y-Y\right\Vert _{\widetilde{W}%
}\left( \left\Vert Y\right\Vert _{\widetilde{W}}+\left\Vert Y_{n}\right\Vert
_{\widetilde{W}}\right) ds| \\
& \leq CM\Vert P_{n}Y-Y\Vert _{L^{2}(\Omega \times (0,T);\widetilde{W}%
)}\left( \Vert Y\Vert _{L^{2}(\Omega \times (0,T);\widetilde{W})}+\Vert
Y_{n}\Vert _{L^{2}(\Omega \times (0,T);\widetilde{W})}\right) \\
& \leq CM\Vert P_{n}Y-Y\Vert _{L^{2}(\Omega \times (0,T);\widetilde{W})},
\end{align*}%
which also converges to zero by (\ref{c02Y}).

Convergences (\ref{c1}) and (\ref{c02Y}) give that
\begin{equation*}
P_{n}Y-Y_{n}\rightarrow 0\qquad \text{ weakly in }L^{2}(\Omega \times (0,T),%
\widetilde{W}),
\end{equation*}%
then for any operator $A\in L^{2}(\Omega \times (0,T),\widetilde{W}^{\ast })$
we have%
\begin{equation*}
\mathbb{E}\int_{0}^{T}\langle A,P_{n}Y-Y_{n}\rangle \,ds\rightarrow 0\qquad
\text{as }n\rightarrow \infty .
\end{equation*}%
Since the function $1_{[0,\tau _{M}]}(s)\xi (s)$ is bounded and independent
of the space variable, we have
\begin{eqnarray*}
&&\Vert 1_{[0,\tau _{M}]}(s)\xi (s)\left( B(Y,Y)-B^{\ast }\right) \Vert
_{L^{2}(\Omega \times (0,T),\widetilde{W}^{\ast })}^{2} \\
&\leq &C\left( \Vert B(Y,Y)\Vert _{L^{2}(\Omega \times (0,T),\widetilde{W}%
^{\ast })}^{2}+\Vert B^{\ast }\Vert _{L^{2}(\Omega \times (0,T),\widetilde{W}%
^{\ast })}^{2}\right) \leq C,
\end{eqnarray*}%
by (\ref{as}), (\ref{l}) and (\ref{c01}). Therefore
\begin{align*}
J_{3}& =2\mathbb{E}\int_{0}^{\tau _{M}}\xi (s)I_{3}\,ds \\
& =2\mathbb{E}\int_{0}^{T}\langle 1_{[0,\tau _{M}]}(s)\xi (s)\left(
B(Y,Y)-B^{\ast }(s)\right) ,P_{n}Y-Y_{n}\rangle \,ds\rightarrow 0\qquad
\text{as }n\rightarrow \infty .
\end{align*}%
We write
\begin{align*}
J_{4}& =\mathbb{E}\int_{0}^{\tau _{M}}\xi (s)\left[ C\left\Vert
P_{n}Y-Y\right\Vert _{V}^{2}+2(P_{n}\widetilde{G}_{n}-P_{n}\widetilde{G}%
^{\ast },P_{n}\widetilde{G}-P_{n}\widetilde{G}^{\ast })_{V}\right] \,ds \\
& =C\mathbb{E}\int_{0}^{T}1_{[0,\tau _{M}]}(s)\xi (s)\left\Vert
P_{n}Y-Y\right\Vert _{V}^{2}ds \\
& +C\mathbb{E}\int_{0}^{T}1_{[0,\tau _{M}]}(s)\xi (s)(\widetilde{G}_{n}-%
\widetilde{G}^{\ast },P_{n}\widetilde{G}-P_{n}\widetilde{G}^{\ast })_{V}\,ds.
\end{align*}%
Due to (\ref{c02Y}), we have
\begin{equation*}
\left\vert \mathbb{E}\int_{0}^{T}1_{[0,\tau _{M}]}(s)\xi (s)\left\Vert
P_{n}Y-Y\right\Vert _{V}^{2}ds\right\vert \leq \mathbb{E}\int_{0}^{T}\left%
\Vert P_{n}Y-Y\right\Vert _{V}^{2}ds\rightarrow 0,\text{ as }n\rightarrow
\infty .
\end{equation*}%
Now, for each stochastic process $Z\in L^{2}(\Omega \times (0,T),H)$ let us
denote by $\widetilde{Z}$ the solution of the modified Stokes problem (\ref%
{GS_NS}). We recall that the operator
\begin{equation*}
A:Z\rightarrow \widetilde{Z}
\end{equation*}%
is linear and continuous operator from $L^{2}(\Omega \times (0,T),V)$ into $%
L^{2}(\Omega \times (0,T),V)$. Applying Proposition A.2 in \cite{B99} (see
also references therein), it follows that $A$ is continuous for the weak
topology, namely if $Z_{n}\rightharpoonup Z$ weakly in $L^{2}(\Omega \times
(0,T),V)$, then $\widetilde{Z}_{n}\rightharpoonup \widetilde{Z}$ weakly in $%
L^{2}(\Omega \times (0,T),V).$ Due to this property and the convergence
result (\ref{c01}), we obtain
\begin{equation}
\widetilde{G}_{n}-\widetilde{G}^{\ast }\rightharpoonup \widetilde{G}-%
\widetilde{G}^{\ast }\text{\qquad\ weakly in }L^{2}(\Omega \times (0,T),V^m).
\label{WCR1}
\end{equation}%
Moreover, we have $\widetilde{G}-\widetilde{G}^{\ast }\in \widetilde{W}^m$ and
\begin{align*}
& P_{n}(\widetilde{G}-\widetilde{G}^{\ast })\rightarrow \widetilde{G}-%
\widetilde{G}^{\ast }\text{\qquad strongly in }L^{2}(\Omega \times (0,T),%
\widetilde{W}^m ), \\
& P_{n}(\widetilde{G}-\widetilde{G}^{\ast })\rightarrow \widetilde{G}-%
\widetilde{G}^{\ast }\text{\qquad strongly in }L^{2}(\Omega \times (0,T),V^m ).
\end{align*}%
Then we can verify that
\begin{equation}
1_{[0,\tau _{M}]}(s)\xi (s)\,P_{n}(\widetilde{G}-\widetilde{G}^{\ast
})\rightarrow 1_{[0,\tau _{M}]}(s)\xi (s)\,(\widetilde{G}-\widetilde{G}%
^{\ast })\text{\qquad strongly in }L^{2}(\Omega \times (0,T),V^m ).
\label{WCR2}
\end{equation}%
As a consequence of (\ref{WCR1}) and (\ref{WCR2}), we have
\begin{align*}
& \mathbb{E}\int_{0}^{T}1_{[0,\tau _{M}]}(s)\xi (s)(\widetilde{G}_{n}-%
\widetilde{G}^{\ast },P_{n}\widetilde{G}-P_{n}\widetilde{G}^{\ast })_{V}\,ds
\\
& =\mathbb{E}\int_{0}^{T}(\widetilde{G}_{n}-\widetilde{G}^{\ast },1_{[0,\tau
_{M}]}(s)\xi (s)P_{n}(\widetilde{G}-\widetilde{G}^{\ast
}))_{V}\,ds\rightarrow 0,\text{ as }n\rightarrow \infty .
\end{align*}%
Collecting all convergence results, we obtain the following strong
convergences\textit{,} depending on the stopping times $\tau _{M},$\
\begin{align*}
\lim_{n\rightarrow \infty }\mathbb{E}\left( \xi (\tau _{M})||P_{n}Y(\tau
_{M})-Y_{n}(\tau _{M})||_{V}^{2}\right) & =0, \\
\lim_{n\rightarrow \infty }\mathbb{E}\int_{0}^{\tau _{M}}\xi (s)\Vert P_{n}%
\widetilde{G}-P_{n}\widetilde{G}^{\ast }\Vert _{V}^{2}ds& =0, \\
\lim_{n\rightarrow \infty }\mathbb{E}\int_{0}^{\tau _{M}}\xi (s)\left( 4\nu
||D(P_{n}Y-Y_{n})||_{2}^{2}+2\nu \gamma \left\Vert P_{n}Y-Y_{n}\right\Vert
_{L^{2}(\Gamma )}^{2}\right) ds& =0,
\end{align*}%
for each $M\in \mathbb{N}$. Since there exists a strictly positive constant $%
\mu,$ such that $\mu\leq 1_{[0,\tau _{M}]}(s)\xi (s)\leq
1, $ it follows that
\begin{equation*}
\lim_{n\rightarrow \infty }\mathbb{E}\int_{0}^{\tau
_{M}}||P_{n}Y-Y_{n}||_{V}^{2}ds=0\quad \text{implying}\quad
\lim_{n\rightarrow \infty }\mathbb{E}\int_{0}^{\tau
_{M}}||Y-Y_{n}||_{V}^{2}ds=0
\end{equation*}%
by (\ref{c02Y}). In addition, considering (\ref{c02}), we have%
\begin{equation}
\mathbb{E}\int_{0}^{\tau _{M}}\Vert \widetilde{G}-\widetilde{G}^{\ast }\Vert
_{V}^{2}ds=0.  \label{RC2}
\end{equation}%
\bigskip

\textit{Step 4. Identification of }
$B^{\ast }(t)$\textit{ with} $B(Y,Y)$\textit{\ and }
$G^{\ast }(t)$ \textit{ with} $G(t,Y).$

Now, we are able to show that the limit function $Y$ satisfies equation (\ref%
{var_form_state}). Integrating equation (\ref{y12}) on the time interval $%
(0,\tau _{M}\wedge t)$, we derive
\begin{eqnarray}
\left( \upsilon \left( Y(\tau _{M}\wedge t)\right) ,\phi \right) -\left(
\upsilon \left( Y_{0}\right) ,\phi \right) &=&\int_{0}^{\tau _{M}\wedge t}
\left[ \left( \nu \Delta Y+U,\phi \right) -\langle B^{\ast }(s),\phi \rangle %
\right] \,ds  \notag \\
&&+\int_{0}^{\tau _{M}\wedge t}\left( G^{\ast }(s),\phi \right) \,dW_{s}
\label{y122}
\end{eqnarray}%
for any $\phi \in \widetilde{W}.$

From (\ref{RC2}) it follows that%
\begin{equation*}
1_{[0,\tau _{M}]}(t)\widetilde{G}=1_{[0,\tau _{M}]}(t)\widetilde{G}^{\ast
}\qquad \text{a. e. in }\Omega \times (0,T),
\end{equation*}%
which implies
\begin{equation}
1_{[0,\tau _{M}]}(t)G(t,Y)=1_{[0,\tau _{M}]}(t)G^{\ast }(t)\qquad \text{a.
e. in }\Omega \times (0,T)  \label{09}
\end{equation}%
by \ (\ref{GS_NS}). Since $%
B(Y_{n},Y_{n})-B(Y,Y)=B(Y_{n},Y_{n}-Y)-B(Y_{n}-Y,Y),$ using (\ref{proB_2})-(%
\ref{proB_3}), we have
\begin{equation*}
\Vert B(Y_{n},Y_{n})-B(Y,Y)\Vert _{\widetilde{W}^{\ast }}\leq C\left( \Vert
Y_{n}\Vert _{\widetilde{W}}+\Vert Y\Vert _{\widetilde{W}}\right) \Vert
Y_{n}-Y\Vert _{V}.
\end{equation*}%
Then for any $\varphi \in L^{\infty }(\Omega \times (0,T);\widetilde{W}),$
using (\ref{as}), (\ref{c1})
\begin{align*}
\bigl\vert\mathbb{E}\int_{0}^{T}1_{[0,\tau _{M}]}(s)\langle B(Y_{n},Y_{n})&
-B(Y,Y),\ \varphi \rangle \,ds\bigr\vert \\
& \leq C\mathbb{E}\int_{0}^{T}1_{[0,\tau _{M}]}(s)\left( \Vert Y_{n}\Vert _{%
\widetilde{W}}+\Vert Y\Vert _{\widetilde{W}}\right) \Vert Y_{n}-Y\Vert
_{V}\Vert \varphi \Vert _{\widetilde{W}}\,ds \\
& \leq C\Vert \varphi \Vert _{L^{\infty }(\Omega \times (0,T),\widetilde{W})}%
\mathbb{E}\int_{0}^{T}\left( \Vert Y_{n}\Vert _{\widetilde{W}}+\Vert Y\Vert
_{\widetilde{W}}\right) \Vert Y_{n}-Y\Vert _{V}\,ds \\
& \leq C\Vert \varphi \Vert _{L^{\infty }(\Omega \times (0,T),\widetilde{W}%
)}\left( \mathbb{E}\int_{0}^{\tau _{M}}\Vert Y_{n}-Y\Vert
_{V}^{2}\,ds\right) ^{\frac{1}{2}}\rightarrow 0,\qquad \text{as }%
n\rightarrow \infty .
\end{align*}%
Taking into account  \ (\ref{c01})$_{1}$ and that the space $L^{\infty }(\Omega \times
(0,T);\widetilde{W})$ is dense in $L^{2}(\Omega \times (0,T);\widetilde{W}),$
we obtain
\begin{equation}
1_{[0,\tau _{M}]}(s)B^{\ast }(s)=1_{[0,\tau _{M}]}(s)B(Y,Y)\qquad \text{a.
e. in }\Omega \times (0,T).  \label{RC2Y}
\end{equation}%
By introducing identities (\ref{09}), (\ref{RC2Y}) in equation (\ref{y122}),
it follows that
\begin{eqnarray}
\left( \upsilon \left( Y(\tau _{M}\wedge t)\right) ,\phi \right) -\left(
\upsilon \left( Y_{0}\right) ,\phi \right)  &=&\int_{0}^{\tau _{M}\wedge t}
\left[ \left( \nu \Delta Y+U,\phi \right) -\langle B(Y,Y),\phi \rangle %
\right] \,ds  \notag \\
&&+\int_{0}^{\tau _{M}\wedge t}\left( G(s,Y),\phi \right) \,dW_{s}.
\label{y1122}
\end{eqnarray}%
Now, reasoning as in (\ref{fc}) \ we have $\tau _{M}\rightarrow T$ \ a. e. in
$\Omega $. We can pass to the limit in each term of equation (\ref{y1122})
in $L^{1}(\Omega \times (0,T))$, as $M\rightarrow \infty $ , by applying the
Lebesgue dominated convergence theorem and the Burkholder-Davis-Gundy
inequality for the last (stochastic) term, deriving equation (\ref%
{var_form_state}) a. e. in $\Omega \times (0,T)$.

Let us notice that the estimates for $Y_n$ in Lemmas \ref{existence_state} and \ref{existence_state1} are valid
also for the limit process $Y$, due to convergence \eqref{c1}.

The uniqueness of the solution $Y$ follows from the stability result that we
will show in the next section. $\hfill \hfill \blacksquare $

\bigskip

\section{Stability result for solutions}

\setcounter{equation}{0}

\bigskip

In this section we will establish a stability property for
solutions of the stochastic second grade fluid model (\ref{equation_etat}).
In spite of the existence result with $H^{3}$ space regularity, the
difference of two solutions can only be estimated (with respect to the
initial data) in space $H^{2}$. It will be convenient to introduce the
following norm on the space $W$
\begin{equation*}
\Vert y\Vert _{W}=\Vert y\Vert _{V}+\Vert \mathbb{P}\upsilon (y)\Vert
_{2},\qquad y\in W.
\end{equation*}%
As a consequence of (\ref{www}) and (\ref{korn}) this norm $\Vert \cdot
\Vert _{W}$ is equivalent to $\Vert \cdot \Vert _{H^{2}}$.

\begin{theorem}
\label{Lips} Assume that for \ some $\ 4\leq p<\infty $
\begin{equation*}
U_{1},U_{2}\in L^{p}(\Omega ,L^{p}(0,T;H(\mathrm{curl};\mathcal{O}))),\qquad
Y_{1,0},Y_{2,0}\in L^{p}(\Omega ,V)\cap L^{2}(\Omega ,\widetilde{W})
\end{equation*}%
and%
\begin{equation*}
Y_{1},Y_{2}\in L^{2}(\Omega ,L^{\infty }(0,T;\widetilde{W}))\cap
L^{p}(\Omega ,L^{\infty }(0,T;V))
\end{equation*}%
are corresponding solutions of (\ref{equation_etat}) in the sense of the
variational equality $(\ref{var_form_state})$.

Then there exist strictly positive constants $C_{3}$ and $C,$\ \ which
depend only on the data (the domain $\mathcal{O},$ \ the regularity of $%
\Gamma ,$ the physical constants $\nu $, $\alpha $, $\gamma $, $K)$ and
satisfy the following estimate
\begin{equation}
\mathbb{E}\sup_{s\in \lbrack 0,t]}\xi (s)\left\Vert
Y_{1}(s)-Y_{2}(s)\right\Vert _{W}^{2}\leq C(\mathbb{E}\left\Vert
Y_{1,0}-Y_{2,0}\right\Vert _{W}^{2}+\mathbb{E}\int_{0}^{t}\xi (s)\left\Vert
U_{1}(s)-U_{2}(s)\right\Vert _{2}^{2}\,ds)  \label{lip_H3}
\end{equation}%
with the function $\xi $ defined as%
\begin{equation*}
\xi (t)=e^{-C_{3}\int_{0}^{t}\left( \left\Vert Y_{1}\right\Vert
_{H^{3}}+\left\Vert Y_{2}\right\Vert _{H^{3}}\right) ds}.
\end{equation*}
\end{theorem}

\textbf{Proof.} The process $Y=Y_{1}-Y_{2}$ satisfies the system
\begin{equation}
\left\{
\begin{array}{ll}
d\upsilon (Y)=(\nu \Delta Y-\mathrm{curl}\,\upsilon (Y)\times Y_{2}-\mathrm{%
curl}\,\upsilon (Y_{1})\times Y-\nabla \pi +U)\,dt &  \\
\qquad \qquad +\left( G(t,Y_{1})-G(t,Y_{2})\right) \,dW_{t} &  \\
\nabla \cdot Y=0 & \quad \mbox{in}\ \mathcal{O}\times (0,T),\vspace{2mm} \\
Y\cdot \mathrm{n}=0,\qquad \left[ 2(\mathrm{n}\cdot DY)+\gamma Y\right]
\cdot \mathrm{\tau }=0 & \quad \mbox{on}\ \Gamma \times (0,T),\vspace{2mm}
\\
Y(0)=Y_{0}=Y_{1,0}-Y_{2,0} & \quad \mbox{in}\ \mathcal{O},%
\end{array}%
\right.  \label{y1-y2}
\end{equation}%
where $\pi =\pi _{1}-\pi _{2}$ and $U=U_{1}-U_{2}$. Applying the operator $%
\left( I-\alpha \mathbb{P}\Delta \right) ^{-1}$ to equation (\ref{y1-y2})$%
_{1}$ we deduce a stochastic differential equation for $Y$, then with the
help of \ It\^{o}'s formula we obtain
\begin{align}
d\left\Vert Y\right\Vert _{V}^{2}& =2\left( \left( \nu \Delta Y-\mathrm{curl}%
\,\upsilon (Y)\times Y_{2}-\mathrm{curl}\,\upsilon (Y_{1})\times Y+U\right)
,Y\right) \,dt\vspace{2mm}  \notag \\
& +\Vert \widetilde{G}_{1}-\widetilde{G}_{2}\Vert _{V}^{2}\,dt+2\left(
G(t,Y_{1})-G(t,Y_{2}),Y\right) \,dW_{t},  \label{EL1}
\end{align}%
where $\widetilde{G}_{i}$ are the solutions of the modified Stokes problem (%
\ref{GS_NS}) with $f=G(t,Y_{i}),\,\ i=1,2.$ Hence, using assumption (\ref{LG}%
), we have
\begin{equation*}
\Vert \widetilde{G}_{1}-\widetilde{G}_{2}\Vert _{V}^{2}\leq C\Vert
G(t,Y_{1})-G(t,Y_{2})\Vert _{2}^{2}\leq C\Vert Y\Vert _{V}^{2}.
\end{equation*}%
Taking into account property (\ref{bb1}), estimate (\ref{rm2_lin}) and the
Young inequality, we derive
\begin{align}
\left\Vert Y(t)\right\Vert _{V}^{2}& +\int_{0}^{t}\left( 4\nu \left\Vert
DY\right\Vert _{2}^{2}+2\nu \gamma \left\Vert Y\right\Vert _{L^{2}(\Gamma
)}^{2}\right) \,ds\leq \left\Vert Y_{0}\right\Vert
_{V}^{2}+C\int_{0}^{t}\left\Vert Y_{2}\right\Vert _{H^{3}}\left\Vert
Y\right\Vert _{Y}^{2}\,ds  \notag \\
& +\int_{0}^{t}\Vert U\Vert _{2}^{2}\,ds+C\int_{0}^{t}\Vert Y\Vert
_{V}^{2}\,ds+2\int_{0}^{t}\left( G(s,Y_{1})-G(s,Y_{2}),Y\right) \,dW_{s}%
\vspace{2mm}.  \label{lip1_V}
\end{align}%
The It\^{o} formula also gives
\begin{align*}
d\left\Vert \mathbb{P}\upsilon \left( Y\right) \right\Vert _{2}^{2}&
=2\left( \nu \Delta Y-\mathrm{curl}\,\upsilon (Y)\times Y_{2}-\mathrm{curl}%
\,\upsilon (Y_{1})\times Y+U,\mathbb{P}\upsilon \left( Y\right) \right) \,dt%
\vspace{2mm} \\
& +\Vert G(t,Y_{1})-G(t,Y_{2})\Vert _{2}^{2}\,dt+2\left(
G(t,Y_{1})-G(t,Y_{2}),\mathbb{P}\upsilon (Y)\right) \,dW_{t}.
\end{align*}%
Estimating the nonlinear term
\begin{equation*}
\left\vert \left( \mathrm{curl}\,\upsilon (Y)\times Y_{2}+\mathrm{curl}%
\,\upsilon (Y_{1})\times Y,\mathbb{P}\upsilon (Y)\right) \right\vert \leq
C_{3}\left( \left\Vert Y_{2}\right\Vert _{H^{3}}+\left\Vert Y_{1}\right\Vert
_{H^{3}}\right) \left( \left\Vert Y\right\Vert _{V}^{2}+\left\Vert \mathbb{P}%
\upsilon (Y)\right\Vert _{2}^{2}\right)
\end{equation*}%
and using (\ref{LG}), we deduce
\begin{align*}
\left\Vert \mathbb{P}\upsilon \left( Y(t)\right) \right\Vert _{2}^{2}+\frac{%
2\nu }{\alpha }\int_{0}^{t}\left\Vert \upsilon \left( Y\right) \right\Vert
_{2}^{2}\,ds& \leq \left\Vert \mathbb{P}\upsilon \left( Y_{0}\right)
\right\Vert _{2}^{2}+2\int_{0}^{t}\left( \frac{\nu }{\alpha }Y+U,\mathbb{P}%
\upsilon \left( Y\right) \right) \,ds \\
& +C_{3}\int_{0}^{t}\left( \left\Vert Y_{2}\right\Vert _{H^{3}}+\left\Vert
Y_{1}\right\Vert _{H^{3}}\right) \left( \left\Vert Y\right\Vert
_{V}^{2}+\left\Vert \mathbb{P}\upsilon (Y)\right\Vert _{2}^{2}\right) \,ds \\
& +C\int_{0}^{t}\Vert Y\Vert _{V}^{2}\,ds+2\int_{0}^{t}\left(
G(s,Y_{1})-G(s,Y_{2}),\mathbb{P}\upsilon (Y)\right) \,dW_{s}.
\end{align*}%
Summing this inequality with (\ref{lip1_V}), we obtain
\begin{align*}
\left\Vert Y(t)\right\Vert _{V}^{2}+\left\Vert \mathbb{P}\upsilon \left(
Y(t)\right) \right\Vert _{2}^{2}& \leq \left\Vert Y_{0}\right\Vert
_{V}^{2}+\left\Vert \mathbb{P}\upsilon \left( Y_{0}\right) \right\Vert
_{2}^{2} \\
& +\int_{0}^{t}\Vert U\Vert _{2}^{2}\,ds+C\int_{0}^{t}\left( \Vert Y\Vert
_{V}^{2}+\Vert \mathbb{P}\upsilon \left( Y\right) \Vert _{2}^{2}\right) \,ds
\\
& +C_{3}\int_{0}^{t}\left( \left\Vert Y_{2}\right\Vert _{H^{3}}+\left\Vert
Y_{1}\right\Vert _{H^{3}}\right) \left( \left\Vert Y\right\Vert
_{V}^{2}+\left\Vert \mathbb{P}\upsilon (Y)\right\Vert _{2}^{2}\right) \,ds \\
& +2\int_{0}^{t}\left( G(s,Y_{1})-G(s,Y_{2}),Y+\mathbb{P}\upsilon (Y)\right)
\,dW_{s}\vspace{2mm}.
\end{align*}%
Taking $\xi (t)=e^{-C_{3}\int_{0}^{t}\left( \left\Vert Y_{2}\right\Vert
_{H^{3}}+\left\Vert Y_{2}\right\Vert _{H^{3}}\right) ds}$ and applying It%
\^{o}'s formula,\ then we easily obtain
\begin{align}
\xi (t)\left( \left\Vert Y(t)\right\Vert _{V}^{2}+\left\Vert \mathbb{P}%
\upsilon \left( Y(t)\right) \right\Vert _{2}^{2}\right) & \leq \left\Vert
Y_{0}\right\Vert _{V}^{2}+\left\Vert \mathbb{P}\upsilon \left( Y_{0}\right)
\right\Vert _{2}^{2}  \notag \\
& +\int_{0}^{t}\xi (s)\Vert U\Vert _{2}^{2}\,ds+C\int_{0}^{t}\xi (s)\left(
\Vert Y\Vert _{V}^{2}+\Vert \mathbb{P}\upsilon \left( Y\right) \Vert
_{2}^{2}\right) \,ds  \notag \\
& +2\int_{0}^{t}\xi (s)\left( G(s,Y_{1})-G(s,Y_{2}),Y+\mathbb{P}\upsilon
(Y)\right) \,dW_{s}.  \label{k}
\end{align}%
The Burkholder-Davis-Gundy inequality gives
\begin{align*}
& \mathbb{E}\sup_{s\in \lbrack 0,t]}\left\vert \int_{0}^{s}\xi (r)\left(
G(r,Y_{1})-G(r,Y_{2}),Y+\mathbb{P}\upsilon (Y)\right) \,dW_{r}\right\vert \\
& \leq \mathbb{E}\left( \int_{0}^{t}\xi ^{2}(s)\left\Vert Y\right\Vert
_{V}^{2}\left\Vert Y+\mathbb{P}\upsilon (Y)\right\Vert _{2}^{2}\,ds\right) ^{%
\frac{1}{2}}\vspace{2mm} \\
& \leq \varepsilon \mathbb{E}\sup_{s\in \lbrack 0,t]}\xi (s)\left(
\left\Vert Y\right\Vert _{V}+\left\Vert \mathbb{P}\upsilon (Y)\right\Vert
_{2}\right) \vspace{2mm} \\
& +C_{\varepsilon }\mathbb{E}\int_{0}^{t}\xi (s)\left( \left\Vert
Y\right\Vert _{V}+\left\Vert \mathbb{P}\upsilon (Y)\right\Vert _{2}\right)
\,ds.
\end{align*}%
Substituting this inequality with $\varepsilon =\frac{1}{2}$ in (\ref{k})
and taking the supremum on the time interval $[0,t]$ and the expectation, we
deduce
\begin{equation*}
\mathbb{E}\sup_{s\in \lbrack 0,t]}\xi (s)\left\Vert Y(s)\right\Vert
_{W}^{2}\leq \mathbb{E}\left\Vert Y_{0}\right\Vert _{W}^{2}+\mathbb{E}%
\int_{0}^{t}\xi (s)\Vert U\Vert _{2}^{2}\,ds+C\mathbb{E}\int_{0}^{t}\xi
(s)\Vert Y(s)\Vert _{W}^{2}\,ds.
\end{equation*}%
Hence Gronwall's inequality yields (\ref{lip_H3}). $\hfill \hfill
\blacksquare $

\bigskip

\bigskip

\textbf{Acknowledgment} We would like to thank the anonymous Referees for
relevant suggestions and comments which contributed to improve the article.

The work of F. Cipriano was partially supported by the Funda\c{c}\~{a}o para
a Ci\^{e}ncia e a Tecnologia (Portuguese Foundation for Science and
Technology) through the project UID/MAT/00297/2013 (Centro de Matem\'{a}tica
e Aplica\c{c}\~{o}es).

\textbf{\bigskip }

\end{document}